\documentclass[12pt]{amsart}

\usepackage{amssymb, stackrel}
\usepackage{mathtools}
\usepackage{txfonts}
\usepackage{eucal}
\usepackage{wasysym}
\usepackage[all]{xy}
\usepackage{relsize}

\usepackage{tikz}

\usepackage{graphicx}

\usepackage{float}

\usepackage{color}

\usepackage{xspace}
\usepackage[all]{xy}
\usepackage[a4paper,body={16.3cm,22.8cm},centering]{geometry}
\usepackage[colorlinks,final,backref=page,hyperindex,pdftex]{hyperref}
\usepackage{enumitem} 

\usetikzlibrary{calc,shapes.geometric}
\tikzset{
baseon/.style={baseline={($(#1)+(0,-0.58ex)$)}},
baseon/.default=current bounding box.center,
every picture/.style=baseon,
lst/.style={},
dst/.style={circle,inner sep=1pt,outer sep=0pt,fill,draw,dst2},
dst2/.style={fill=white},
ddst/.style={diamond,draw,inner sep=1pt},
eest/.style={ellipse,draw,inner sep=1pt,minimum size=2ex},
}

\font \sevenrm=cmr7
\font \fiverm=cmr5
\newcommand{\nc}{\newcommand}
\nc\smsc{0.8}

\nc{\ooverline}[1]{\overline{\overline #1}}
\nc{\uunderline}[1]{\underline{\underline #1}}
\def \restr#1{\mathstrut_{ | \, \scriptstyle{ #1}}}
\def \srestr#1{\mathstrut_{\scriptstyle |}\hbox to
-1.5pt{}\raise-4pt\hbox{$\hskip 1pt\scriptscriptstyle #1$}}
\nc{\mop}[1]{\mathop{\hbox {\rm #1} }\nolimits}
\nc{\gmop}[1]{\mathop{\hbox {\bf #1} }\nolimits}

\nc{\smop}[1]{\mathop{\hbox {\sevenrm #1} }\nolimits}
\nc{\ssmop}[1]{\mathop{\hbox {\fiverm #1} }\nolimits}
\nc{\mopl}[1]{\mathop{\hbox {\rm #1} }\limits}
\def\dbar{d\hskip-3pt \raise 4pt\hbox{-}}
\nc{\smopl}[1]{\mathop{\hbox {\sevenrm #1} }\limits}
\nc{\ssmopl}[1]{\mathop{\hbox {\fiverm #1} }\limits}
\newcommand{\delete}[1]{}

\nc{\mlabel}[1]{\label{#1}}  
\nc{\mcite}[1]{\cite{#1}}  
\nc{\mref}[1]{\ref{#1}}  
\nc{\mbibitem}[1]{\bibitem{#1}} 

\delete{
\nc{\mlabel}[1]{\label{#1}  
{\hfill \hspace{1cm}{\small\tt{{\ }\hfill(#1)}}}}
\nc{\mcite}[1]{\cite{#1}{\small{\tt{{\ }(#1)}}}}  
\nc{\mref}[1]{\ref{#1}{{\tt{{\ }(#1)}}}}  
\nc{\mbibitem}[1]{\bibitem[\bf #1]{#1}} 
}

\newtheorem{theorem}{Theorem}[section]

\newtheorem{notation}[theorem]{Notation}
\newtheorem{corollary}[theorem]{Corollary}
\newtheorem{proposition}[theorem]{Proposition}
\newtheorem{lemma}[theorem]{Lemma}

\theoremstyle{remark}
\newtheorem{remark}[theorem]{Remark}
\newtheorem{example}[theorem]{Example}
\numberwithin{equation}{section}

{\theoremstyle{definition}
\newtheorem{definition}[theorem]{Definition}
}

\newcommand\alphlist{a,b,c,d,e,f,g,h,i,j,k,l,m,n,o,p,q,r,s,t,u,v,w,x,y,z}
\newcommand\Alphlist{A,B,C,D,E,F,G,H,I,J,K,L,M,N,O,P,Q,R,S,T,U,V,W,X,Y,Z}
\newcommand\getcmds[3]{\expandafter\newcommand\csname #2#1\endcsname{#3{#1}}}
\makeatletter
\@for\x:=\alphlist\do{\expandafter\getcmds\expandafter{\x}{frak}{\mathfrak}}
\@for\x:=\Alphlist\do{\expandafter\getcmds\expandafter{\x}{frak}{\mathfrak}}
\makeatother

\nc{\bfk}{{\bf k}}
\nc{\sha}{\shuffle}

\nc{\id}{\mathrm{id}}
\nc{\Id}{\mathrm{Id}}
\nc{\lbar}[1]{\overline{#1}}
\nc{\ot}{\otimes}
\nc{\dep}{\mathrm{dep}}
\nc{\ver}{\mathrm{ver}}

\nc{\tred}[1]{\textcolor{red}{#1}} \nc{\tgreen}[1]{\textcolor{green}{#1}}
\nc{\tblue}[1]{\textcolor{blue}{#1}} \nc{\tpurple}[1]{\textcolor{purple}{#1}}
\nc{\tcyan}[1]{\textcolor{cyan}{#1}} 
\nc{\tblk}[1]{\textcolor{black}{#1}}

\nc{\xadd}[1]{\tred{#1}}
\nc{\xing}[1]{\tblue{\underline{Xing:}#1 }}
\nc{\nan}[1]{\tblue{\underline{Nan:}#1 }}
\nc{\dominique}[1]{\tpurple{\underline{Dominique: }#1 }}
\long\def\ignore#1{}
\nc{\zhicheng}[1]{\tpurple{\underline{Zhicheng: }#1 }}
\nc{\dz}[1]{\tpurple{\underline{Dominique and Zhicheng: }#1 }}
\long\def\ignore#1{}

\def\zzz#1`#2...#3`#4...#5`#6@{%
--++(#1)
node[dst,label={#5:$#6$},name=#2]{}
node[midway,auto,#3]{$#4$}
}
\def\ddd#1`#2`#3@{+(#1)node[ddst,name=#2]{$#3$}}
\def\eee#1`#2`#3@{+(#1)node[eest,name=#2]{$#3$}}
\def\xxx#1`#2@{node[midway,auto,inner sep=1pt,#1]{$#2$}}
\def\pp#1`#2`#3@{node[dst,label={#2:$#3$},pos=#1]{}}
\def\oo#1`#2`#3@{\path (o) node[dst,label={#2:$#3$},name=o,#1]{};}
\def\eoo#1`#2@{\node[eest,name=o,#1] at (o) {$#2$};}

\newif\ifshowjdq
\showjdqtrue
\newcommand\setXXclip[3]{%
\def\XXheight{#1}\def\XXdepth{#2}\def\XXwidth{#3}}
\setXXclip{1}{-0.5}{1.1}

\makeatletter
\newcommand\simra{\mathrel{\mathpalette\@verra\sim}}
\def\@verra#1#2{\lower.5\p@\vbox{\lineskiplimit\maxdimen \lineskip-.5\p@
\ialign{$\m@th#1\hfil##\hfil$\crcr#2\crcr\rightarrow\crcr}}}
\makeatother

\nc{\dnx}{\Delta_n A} \nc{\dx}{\Delta A} \nc{\dgp}{{\rm deg_{P}}}
\nc{\dgt}{{\rm deg_{T}}} \nc{\dg}{{\rm deg}} \nc{\ida}{ID($A$)} \nc{\tu}{\tilde{u}} \nc{\tv}{\tilde{v}}
\nc{\nr}{\calr_n} \nc{\nz}{\calz_n} \nc{\fun}{\cala_{n,d}}
 \nc{\fbase}{\calb} \nc{\LF}{\mathrm{RF}} \nc{\FFA}{\mathrm{LF}} \nc{\irr}{\mathrm{Irr}}
 \nc{\result}{\bfk\mathrm{Irr}(S_n)}  \nc{\I}{I_{\mathrm{ID},n}^0}
 \nc{\nrs}{\calr_n^\star} \nc{\ii}{\mathrm{I}} \nc{\iii}{\mathrm{II}}
\nc{\intl}{{\rm int}}\nc{\ws}[1]{{#1}}\nc{\deleted}[1]{\delete{#1}}\nc{\plas}{placements\xspace}

\nc{\bim}[1]{#1}  \nc{\shaop}{\sha_{\Omega}^{+}}  \nc{\shao}{\sha_{\Omega}}
\nc{\bbim}[2]{#1 #2} \nc{\bbbim}[2]{#1,\, #2} \nc{\RBF}{{\rm RBF}}
\nc{\frb}{F_{\RB}} \nc{\shaf}{\ssha_{\tiny{\Omega}}} \nc{\sham}{\diamond_{\tiny{\Omega}}}
\nc{\lf}{\lfloor} \nc{\rf}{\rfloor} \nc{\shan}{\ssha_{\lambda}}
\nc{\rlex}{{\rm {lex}}} \nc{\bb}{\Box} \nc{\ra}{\rightarrow}
\nc{\e}{{\rm {e}}}
\nc{\DDF}{\mathrm{DD}(X,\,\Omega)}\nc{\DTF}{\mathrm{DT}(X,\,\Omega)} \nc{\DT}{\mathrm{DT}'(\Omega,\,V)}
\nc{\bra}{\mathrm{bra}} \nc{\bre}{\mathrm{bre}}
\nc{\dec}{\mathrm{dec}} \nc{\diamondw}{\diamond_{w}}
\nc{\type}{\mathrm{type}}

\nc\caF[1]{\cal{F}_{#1}(X,\,\Omega)}
\nc\calt{\cal{T}(X,\,\Omega)} \nc\caltn{\cal{T}_n(X,\,\Omega)}
\nc\calta{\cal{T}_0(X,\,\Omega)}
\nc\caltb{\cal{T}_1(X,\,\Omega)}
\nc\caltc{\cal{T}_2(X,\,\Omega)}
\nc\caltd{\cal{T}_3(X,\,\Omega)}
\nc\caltm{\cal{T}_m(X,\,\Omega)}
\nc\caltx{\cal{T}(X)}
\nc\calf{\cal{F}(X,\,\Omega)}
\nc\fram{\frak{M}(\Omega,\, X)}
\nc\shaw{\sha^{NC}_w(\Omega,\, X)}
\nc\dw{\diamond_w} \nc\dl{\diamond_\ell}
\nc\shal{\sha^{NC}_\ell(X,\, \Omega)} \nc\shav{\sha^{NC}_w(\Omega,\, V)} \nc\shat{\sha^{NC,1}_w(\Omega,\, T^{+}(V))}
\nc{\cfo}{\cal{F}(X,\,\Omega)}
\nc{\sh}{\rm{Sh}}
\nc{\lar}{\varinjlim}
\nc\XO{(X,\,\Omega)}
\def\cxo#1#2;{\cal{#1}#2\XO}
\nc\lrf[2]{B_{#2}^+(#1)}
\nc{\fd}{\mathrm{\text{typed angularly decorated planar rooted trees}}}
\nc{\rb}{\mathrm{RBFWs}} \nc{\dfw}{\mathrm{DFW{(X)}}} \nc{\tfw}{\mathrm{TFW{(X)}}}
\nc{\tfv}{\mathrm{TFW{(V)}}}

\def\Ve#1,#2,#3;{\vee_{#1,\,(#2,\,#3)}}
\def\bigv#1;#2;#3;{\bigvee\nolimits_{#1}^{#2;\,#3}}
\nc\rjt[2]{\mathrel{\mathop{\longrightarrow}\limits^{#1\hfill}_{\hfill#2}}}
\nc{\pl}{\cal{PLF}}
\nc{\tr}{\cal{RTF}}
\nc{\im}{\mathrm{Im} \, }
\nc{\ff}{\cal{F}_\Omega}
\nc{\tm}{T_\Omega}
\nc{\calp}{\cal{P}}
\makeatletter
\nc\dd{\@ifnextchar'{\ddA}{\ddB}}
\def\ddA'#1;{\rhd'_{#1\,}}
\def\ddB#1;{\rhd_{#1\,}}
\nc{\pbt}{\mathrm{PBT}}
\nc{\ad}{\mathrm{ad}}

\nc{\hRR}{\hat {\mathbb R}}
\nc{\RR}{{\mathbb R}}
\nc{\NN}{{\mathbb N}}
\nc{\ol}[1]{\overline{#1}}
\nc{\set}{{\rm \bf set}}
\nc{\setx}{{\rm \bf set^{\times}}}
\nc{\setn}{{\rm \bf set_{{\mathbb N}}}}
\nc{\Sb}{{\mathbb S}}
\nc{\cp}{\Vdash}
\nc{\te}{\otimes}
\nc{\sus}{\subseteq}
\nc{\BF}{{\rm BF}\,}
\nc{\Pio}{\Pi^\circ}
\nc{\SM}{{\rm SM}\,}
\nc{\ti}{\times}
\nc{\ben}{{\bf 1}}
\nc{\bigsum}{{\mathlarger{\sum}}}
\nc{\nil}{{\mathbf 0}}
\nc{\Hom}{\text{Hom}}
\nc{\cT}{{\mathcal T}}

\definecolor{tol1}{HTML}{332288}
\definecolor{tol2}{HTML}{AA4499}
\definecolor{tol3}{HTML}{DDCC77}
\definecolor{tol4}{HTML}{44AA99}
\definecolor{tol5}{HTML}{88CCEE}
\definecolor{tol6}{HTML}{CC6677}
\definecolor{tol7}{HTML}{117733}
\definecolor{tol8}{HTML}{882255}

\renewcommand{\cT}{{\mathcal T}}

\newcommand{\rleftarrows}{\mathrel{\raise.75ex\hbox{\oalign{%
  $\hfil\scriptstyle\relbar$\cr
  \vrule width0pt height.5ex$\scriptstyle\smash\leftarrow$\cr}}}}
\newcommand{\rightlarrows}{\mathrel{\raise.75ex\hbox{\oalign{%
  $\scriptstyle\rightarrow$\hfil\cr
  $\scriptstyle\vrule width0pt height.5ex\relbar$\cr}}}}
\newcommand{\Rrelbar}{\mathrel{\raise.75ex\hbox{\oalign{%
  $\scriptstyle\relbar$\cr
  \vrule width0pt height.5ex$\scriptstyle\relbar$}}}}

\def\shu{{\scriptstyle \sqcup\hskip -1.4pt\sqcup}}
\newcommand{\calH}{\mathcal{H}}

\makeatother
\begin{document}
\title[Controlled rough paths: a general Hopf-algebraic setting]
{Controlled rough paths: a general Hopf-algebraic setting}
\thispagestyle{empty}

\author{Zhicheng Zhu, Xing Gao, Nannan Li,  Dominique Manchon}
\address{School of Mathematics and Statistics, Lanzhou University,
Lanzhou, 730000, China}
\email{zhuzhch16@lzu.edu.cn}
\email{linann0601@163.com}

\address{School of Mathematics and Statistics, Lanzhou University,
Lanzhou, 730000, China; Gansu Provincial Research Center for Basic Disciplines of Mathematics and Statistics, Lanzhou, 730070, China}
\email{gaoxing@lzu.edu.cn}

\address{Laboratoire de Math\'ematiques Blaise Pascal,
CNRS--Universit\'e Clermont-Auvergne,
3 place Vasar\'ely, CS 60026,
F63178 Aubi\`ere, France}
\email{Dominique.Manchon@uca.fr}
\tikzset{
stdNode/.style={rounded corners, draw, align=right},
greenRed/.style={stdNode, top color=green, bottom color=red},
blueRed/.style={stdNode, top color=blue, bottom color=red}
}
	
\date{\today}
\begin{abstract}
We set up controlled rough paths for a class of combinatorial Hopf algebras, encompassing shuffle, Butcher-Connes-Kreimer and Munthe-Kaas--Wright Hopf algebras. The class of controls we consider encompasses both H\"older continuous paths and (not necessarily continuous) paths with bounded $p$-variation.  We prove existence and uniqueness of the solution of a lifted initial value problem in this general setting by applying the fixed point method in a suitable Banach space of controlled rough paths, and we prove a universal limit theorem addressing the robustness of the solution with respect to the parameters and the initial condition.
\end{abstract}

\makeatletter
\@namedef{subjclassname@2020}{\textup{2020} Mathematics Subject Classification}
\makeatother
\subjclass[2020]{
60L20, 
60L50, 
60H99, 
34K50, 
37H10,  
16T05.  	
}

\keywords{Truncated Hopf algebra, truncated (controlled) rough path, rough integral, universal limit theorem.}
\maketitle

\setcounter{tocdepth}{1}
\tableofcontents



\section{Introduction}
In this paper, the universal limit theorem for controlled rough paths is generalized to driving rough paths relative to a broad class of commutative connected graded Hopf algebras. This framework encompasses weakly geometric, branched and planarly branched rough paths.
\subsection{Rough paths and branched rough paths}
Rough path theory, introduced by T. Lyons in 1998~\cite{L98}, provides a robust mathematical framework for analyzing rough differential equations (RDEs) driven by highly irregular signals, beyond the reach of classical stochastic or Young integration. These equations take the form:
\begin{equation} \mlabel{ivp0}
dZ_{t} = \sum\limits_{i=1}^{d} \varphi_i(Z_{t}) dX_{t}^{i}, \quad \forall t \in [0, T],
\end{equation}
where $X:[0,T]\to \mathbb R^d$ is a path that may exhibit low regularity, and $(\varphi_i)_{i=1,\ldots, d}$ is a collection of smooth vector fields on $\mathbb R^e$. These RDEs are typically interpreted in the integral form:
\begin{equation}\mlabel{ivp}
 Z_t = Z_0 + \int_0^t \sum_{i=1}^d \varphi_i(Z_s) dX^i_s, \quad \forall t \in [0, T],
\end{equation}
highlighting the need to define integrals against possibly non-smooth paths $X$.
The innovation of rough path theory lies in its algebraic enhancement of the path
$X$ by encoding information about its iterated integrals---essentially, the higher-order interactions needed to define and solve RDEs~(\mref{ivp0}). When $X$ is only $\gamma$-H\"older continuous for some $\gamma\in ]0,1]$,
or more generally of bounded $1/\gamma$-variation, the classical notion of iterated integrals breaks down. To address this, Lyons' rough path theory introduces algebraic substitutes for Chen's iterated integrals that remain meaningful in this rough setting.\\

Informally, a rough path $\mathbf X$ above a path $X$ consists of a two-parameter family $(\mathbf X_{st})_{s,t\in [0,T]}$
of linear functionals on the (suitably truncated) shuffle Hopf algebra generated by words over the alphabet $[d]:=\{1,\ldots, d\}$.
These functionals are designed so that evaluating $\mathbf X_{st}$ on words of length one recovers the increments of the path $X$. Furthermore, they satisfy Chen's lemma
$\mathbf X_{st}=\mathbf X_{su}\star\mathbf X_{ut}$, along with appropriate analytic estimates.
When further each $\mathbf X_{st}$ acts as a truncated character---meaning that for any pair $(x,y)$ of words, the relation
$\mathbf X_{st}(x\shu y)=\mathbf X_{st}(x)\mathbf X_{st}(y)$
is valid whenever the sum of their length does not exceed the truncation parameter---the rough path $\mathbf X$ is called weakly geometric~\cite{FH2,Geng, Kel} (for the subtle difference between \textsl{weakly geometric} and \textsl{geometric}, see \cite[Remark 1.3]{HK2015}).
A noteworthy milestone in rough path theory is M. Hairer's use~\cite{Hair13a} of rough paths to solve the Kardar-Parisi-Zhang equation~\cite{C2012,KPZ}. See also~\mcite{Hair11, Hair13b, LG25} for the application of rough paths in solving Burgers-like equations.\\

In contrast to what happens for a differentiable path---whose signature (the collection of its Chen's iterated integrals) is uniquely defined---the construction of rough paths above a given $\gamma$-H\"older continuous path is not unique when $\gamma <1/2$~\cite{LV07, TZ20}. Nevertheless, the sewing lemma~\cite{FP2006} provides a powerful structural result: it guarantees that it suffices to define $\mathbf X$
on words of length at most $N=\lfloor1/\gamma\rfloor$;
the values on longer words are then determined recursively and uniquely.
This principle also extends beyond continuous paths. It applies, for instance, to paths with jumps that possess bounded $1/\gamma$-variation, opening the door to a broader class of driving paths $X$ in RDEs~(\mref{ivp0}), including those arising in models with discontinuities or hybrid dynamics.\\

Branched rough paths, introduced by M. Gubinelli~\cite{Gub1}, extend the concept of weakly geometric rough paths by replacing the algebraic structure underlying their definition. Instead of working with the shuffle Hopf algebra of words formed from the alphabet $[d] := {1, \ldots, d}$, branched rough paths are built upon the Butcher-Connes-Kreimer Hopf algebra of rooted forests, where each node is decorated by an element of $[d]$. This shift from linear words to tree-like structures allows for a more flexible and natural representation of certain algebraic and combinatorial features arising in rough integration.
Given a rough path, one can canonically construct a corresponding branched rough path by precomposing the linear functionals $\mathbf{X}_{st}$ with the arborification map $\mathfrak{a}$ restricted to trees, and extending multiplicatively. This arborification map $\mathfrak{a}$ is a surjective morphism of Hopf algebras from the Butcher-Connes-Kreimer algebra onto the shuffle Hopf algebra, encoding the process of flattening trees into words while preserving algebraic structure.
Importantly, any branched rough path can also be interpreted as a weakly geometric rough path by lifting it to the shuffle algebra over an infinite alphabet~\cite{Fer1, Kel}.\\

Planarly branched rough paths are defined in close analogy to branched rough paths, with the crucial distinction that the underlying algebraic structure is given by the Munthe-Kaas--Wright Hopf algebra of $[d]$-decorated planar rooted forests, rather than the Butcher-Connes-Kreimer Hopf algebra. In this setting, the planarity of the forests introduces an additional layer of combinatorial structure, as the ordering of branches now carries meaningful information. This richer algebraic framework captures more refined dependencies and is particularly well-suited for encoding operations that are sensitive to the order of compositions.
Planarly branched rough paths have proven especially effective in the context of RDEs~(\mref{ivp0}) on homogeneous spaces, where the symmetry and geometry of the underlying space demand a more structured approach~\mcite{CFMK}.
Moreover, the theory of regularity structures~\cite{BH, BHZ, CW17, Hair14, LOTT}, initiated by M. Hairer to handle singular stochastic partial differential equations, has been extended to incorporate planar rooted trees and forests. This adaptation, discussed in~\cite{R22}, enables the use of planarly branched structures within the broader framework of planar regularity structures, further highlighting their versatility and potential in analyzing complex, highly irregular systems.

\subsection{Controlled (branched) rough paths and universal limit theorem}
Controlled rough paths, a powerful concept introduced by M. Gubinelli, were first developed in the shuffle algebra setting~\cite{Gub04} and later extended to the branched setting~\cite{Gub1} (see also~\cite{Kel,MSAN} for further developments). Given a weakly geometric (resp. branched, resp. planarly branched) rough path $\mathbf X$ lying above a path $X$, an $\mathbf X$-controlled rough path above $Z=(Z^1,\ldots,Z^e):[0,T]\to \mathbb R^e$ is defined as a collection $\mathbf Z=(\mathbf Z^1,\ldots ,\mathbf Z^e)$, where each $\mathbf Z^i$ is a map from $[0,T]$
into a (suitably truncated) Hopf algebra of words (resp. rooted forests, resp. planar rooted forests), such that
\[\langle h, \,\mathbf Z_t^i\rangle=\langle \mathbf X_{st}\star h,\,\mathbf Z^i_s\rangle + \hbox{ \sl small \rm remainder},\]
and the path $Z^i$ is recovered via the counit map $\epsilon$ through the identity
$Z^i=\epsilon\circ \mathbf Z^i$.
The space of $\mathbf X$-controlled rough paths on $[0,T]$ forms a Banach space, which is stable under smooth functions:  any smooth map $\varphi:\mathbb R^e\to \mathbb R^e$ gives rise to a unique controlled rough path $\varphi(\mathbf Z)$ above the path $\varphi\circ Z$.
This framework enables us to lift~\eqref{ivp} into the space of controlled rough paths, yielding the lifted integral equation:
\begin{equation}\label{ivplifted}
\mathbf Z_t=\mathbf Z_0+\sum_{i=1}^d\int_0^t  \varphi_i(\mathbf Z)_s\, d X^i_s.
\end{equation}
Here, the integrals $\int_0^t\varphi_i(\mathbf Z)_s\,dX_s^i$  are defined in a precise and algebraic way for any
$\mathbf X$-controlled rough path $\mathbf Z$; see~\cite[(3.25)]{Kel} for full details.\\

As explained in~\cite[Paragraph 3.3.7]{Kel}, the Banach fixed point theorem, applied to the
transformation
\begin{equation}
\mathcal M_T:\mathbf Z \longmapsto \mathbf Z_0 +\sum_{i=1}^d \int_0^\bullet\varphi_i(\mathbf Z)_s  \, dX^i_s
\mlabel{eq:firstm}
\end{equation}
of the Banach space of $\mathbf X$-controlled rough paths on $[0,T]$,
yields existence and uniqueness for solutions of the lifted integral equation~\eqref{ivplifted}, provided $T$ is small enough. The {\em universal limit theorem} ensures existence and uniqueness of the solution even for a bigger bound $T$, as well as its robustness with respect its parameters: initial condition, vector fields $\varphi_i$ and driving rough path $\mathbf X$. In greater detail, the universal limit theorem was originally established by T. Lyons~\cite[Theorem 4.1.1]{L98} weakly geometric rough paths of roughness $\frac{1}{3}< \gamma \leq \frac{1}{2}$, where it emerged as a foundational result in the theory of rough paths. The terminology was coined by P. Malliavin~\cite{M97}, who recognized its significance in stochastic analysis. See also~\cite{FV2,Geng} for details.
Subsequently, in the same year, Lyons and Qian~\cite{LQ98} refined and generalized the universal limit theorem to encompass all Lyons' rough paths with roughness $\frac{1}{3}< \gamma \leq \frac{1}{2}$, not just the weakly geometric case. This development widened the applicability of the universal limit theorem significantly.
Later, in~\cite{LQ02}, Lyons and Z. Qian extended the result even further by proving the universal limit theorem for geometric rough paths with arbitrary roughness $0< \gamma \leq 1$, thus covering the full spectrum of H\"older-like regularities encountered in practice.
More recently, Friz and Zhang~\cite{FZ18} adapted the universal limit theorem to the branched rough path setting, proving its validity for any roughness $0< \gamma \leq 1$ in the bounded $1/\gamma$-variation framework, which notably accommodates paths with jumps---a key advancement in the theory's applicability to non-smooth dynamics.
Finally, the second, third, and fourth authors of the present paper extended the theorem to the planarly branched rough path framework for H\"older exponents in the lower roughness regime $1/4<\gamma\le 1/3$~\cite{GLM24}.
This further broadens the scope of the rough path machinery to more complex combinatorial and geometric structures.\\

It is noteworthy that the foundational idea underlying the universal limit theorem had already found application in the theory of stochastic differential equations (SDEs) prior to the formal development of rough path theory in~\cite{L98}. In fact, both It$\mathrm{\hat{o}}$-type and Stratonovich-type SDEs can be regarded as particular instances of the RDEs~\eqref{ivp0}, when driven by stochastic processes such as Brownian motion~\mcite{FGL15}. A pioneering contribution in this direction was made by E. Wong and M. Zakai in their seminal works~\cite{WZ1, WZ2}, where they established an early form of a universal limit theorem for Stratonovich-type SDEs. Their results demonstrated that solutions to SDEs driven by smooth approximations of Brownian motion converge---under suitable conditions---to the Stratonovich solution of the corresponding stochastic equation. This insight not only provided a rigorous justification for using Stratonovich integrals in physical models but also foreshadowed the algebraic and analytic structure later formalized in the rough path framework.
Further developments and discussions of these ideas can be found in~\cite{D07, M97, Pri}, which highlight the historical and conceptual continuity between stochastic analysis and the modern theory of rough paths. In this light, rough path theory may be seen as a natural extension of classical stochastic calculus, offering a unified and deterministic approach to integration along irregular paths, whether arising from randomness or other forms of roughness.

\subsection{Motivation and outline of the paper}
Our motivation is to describe an algebraic framework, as broad as possible, in which a proof of the universal limit theorem using the Banach fixed point theorem can be implemented. The analytical framework, based on controls fulfilling a right-continuity property, contains both H\"older-continuous paths and paths with bounded $p$-variation, possibly with jumps. Our approach encompasses in particular the universal limit theorem for branched rough paths in the latter case \cite[Theorems 4.19 and 4.20]{FZ18}.
\vspace{0.4cm}

The present paper is organized as follows. We define truncated connected graded Hopf algebras in Section \ref{sect:truncated}, and prove that it is a self-dual notion (Remark \ref{self-dual}). We recall from \cite[Paragraph 1.1]{FZ18} the definition and the main properties of controls in Section \ref{sect:controls}, and introduce our family of right-continuous controls in Definition \ref{rc-controls}. In Section \ref{sect:crp}, we first introduce truncated rough paths associated to a truncated connected graded Hopf algebra, and we define controlled rough paths driven by such a truncated rough path in Paragraph \ref{ssect:trp}. Banach spaces of controlled rough paths are introduced  in the spirit of \cite[Paragraph 4.3]{FZ18}. We then define the transformation of a controlled rough path by a smooth map $\varphi$ in Paragraph \ref{ssect:phi}, and prove the continuity of this transformation (Proposition \ref{main-est-varphi}). We end up Section \ref{sect:crp} with Paragraph \ref{ssect:integration} devoted to the sewing lemma in the general setting of \cite{FZ18} (Lemma 2.5 therein), and to integration of controlled rough paths.\\

Section \ref{sect:ult} first establishes the local existence and uniqueness of solutions to the initial value problem \eqref{ivplifted} in a suitable Banach space of controlled rough paths (Theorem \ref{thm:local}). This is achieved via the Banach fixed point theorem, enabled by the contraction Estimate \eqref{contraction}. Building on this result, we then derive our main theorem in this section (Theorem \ref{thm:ult}), which asserts the global existence and uniqueness of solutions to the initial value problem \eqref{ivplifted}. The last Section \ref{ult-bis} addresses the universal limit theorem (Theorem \ref{thm:ultfinal}), i.e. the continuity of the solution of \eqref{ivplifted} with respect to the initial value, the vector fields $(\varphi_i)_{i=1,\ldots, d}$, and the driving rough path $\mathbf X$.\\

\noindent{\bf Acknowledgements:} We warmly thank Xi Geng for having pointed out Reference \cite{FZ18} to our attention, as well as Yingtong Hou and Huilin Zhang for illuminating discussions. The second author is supported by the National Natural Science Foundation of China (12571019), the Natural Science Foundation of Gansu Province (25JRRA644), Innovative Fundamental Research Group Project of Gansu Province (23JRRA684) and Longyuan Young Talents of Gansu Province.

\begin{notation}\label{notation-intro}\rm
\ignore{Denote by $\RR$ the field of real numbers and by $M( \bullet )$ a universal function that is increasing in all variables. The function of $M(\bullet) $ may differ from line to line.} Throughout this paper, we fix two positive integers $d$ and $e$. The first integer is the number of vector fields driving~\eqref{ivp0}, the second integer is the dimension of the $\RR$-vector space in which the solutions of \eqref{ivp0} take their values.
\end{notation}
\section{Preliminaries on normed vector spaces}\label{sect:nvs}
Let $V$ and $W$ be two normed vector spaces, and let $A:V\to W$ be a linear map. The operator norm of $A$ is defined by
\[\|A\|_{V\to W}:=\mopl{sup}_{x\in V,\, \|x\|\neq 0}\frac{\|A(x)\|}{\|x\|}.\]
Equivalently,
\[\|A\|_{V\to W}:=\mopl{sup}_{x\in V,\, \|x\|=1}\|A(x)\|.\]
If $X$ is another normed vector space and $B:W\to X$ another linear map, the multiplicative property
\[\|B\circ A\|_{V\to X}\le \|A\|_{V\to W}\|B\|_{W\to X}\]
holds. For any normed vector space $V$, let us denote by $V^*$ (resp. $V'$) be its algebraic (resp. topological) dual, i.e. the space of linear (resp. bounded linear) maps form $V$ to $\mathbb R$. Here $\mathbb R$ is considered as a one-dimensional normed vector space, with $\|t\|=|t|$ for any real $t$. The topological dual is endowed with the operator norm $\xi\mapsto \|\xi\|_{V\to\mathbb R}$.
\begin{proposition}\label{norm-transpose}
Let $A:V\to W$ be a bounded linear map. Its transpose $A^*:W^*\to V^*$ restricts to a bounded linear map $A':W'\to V'$ and the equality
\[\|A'\|_{W'\to V'}=\|A\|_{V\to W}\]
holds.
\end{proposition}

\begin{proof}
The following proof of this well-known fact is borrowed from \cite{AB1999}
\begin{eqnarray*}
\|A'\|_{W'\to V'}&=&\mopl{sup}_{\xi\in W',\, \|\xi\|= 1}\|A'\xi\|\\
&=&\mopl{sup}_{\xi\in W',\, \|\xi\|= 1}\, \mopl{sup}_{x\in V,\, \|x\|= 1}|\langle A'\xi,\, x\rangle|\\
&=&\mopl{sup}_{x\in V,\, \|x\|= 1}\, \mopl{sup}_{\xi\in W',\, \|\xi\|= 1}|\langle \xi,\, Ax\rangle|\\
&=&\mopl{sup}_{x\in V,\, \|x\|= 1}\|Ax\|\hspace{2cm} \hbox{ (from the Hahn-Banach theorem)}\\
&=&\|A\|_{V\to W}.  
\end{eqnarray*}
\end{proof}

\noindent If $V$ and $W$ are finite-dimensional, we of course have $V'=V^*$, $W'=W^*$ and $A'=A^*$.
\begin{definition}
Let $\|-\|_V$ be a norm on a finite-dimensional vector space $V$, let $\|-\|_{V^*}$ be the corresponding operator norm on $V^*$, and let $\|-\|_{V\otimes V}$ be the operator norm $\|-\|_{V^*\to V}$, where $V\otimes V$ is naturally identified with the space of linear maps from $V^*$ to $V$. The \textsl{norm} on $V\otimes V$ corresponding to $\|-\|_V$ is defined by\footnote{The case of infinite-dimensional normed vector spaces is much more subtle. See e.g. \cite{G1954}.}
\begin{equation}\label{norm-tens}
\|A\|_{V\otimes V} :=\mopl{sup}_{\|\xi\|_{V^*}=1}\|A\xi\|_V=\mopl{sup}_{\|\xi\|_{V^*}=1}\,\mopl{sup}_{\|\eta\|_{V^*}=1}|\langle A\xi,\eta\rangle|.
\end{equation}
In particular, for any $u,v\in V$ we have
\[\|u\otimes v\|_{V\otimes V}=\|u\|_V \|v\|_V.\]
\end{definition}
\begin{corollary}\label{good-norm}
Let $V$ be a finite-dimensional vector space, and let $\Delta:V\to V\otimes V$ be a nonzero linear map. There exists a norm $\|-\|_V$ on $V$ such that $\|\Delta\|_{V\to V\otimes V}= 1$. Let
\begin{eqnarray*}
m:=\Delta^*:V^*\otimes V^*&\longrightarrow & V^*\\
\xi\otimes\eta &\longmapsto& \xi\star\eta
\end{eqnarray*}
be the transpose of $\Delta$. Then, for this choice of norm on $V$ we have
\[\|\xi\star\eta\|_{V^*}\le \|\xi\|_{V^*}\|\eta\|_{V^*}\]
for any $\xi,\eta\in V^*$.
\end{corollary}
\begin{proof}
Let $N_V$ be an arbitrary norm on $V$, let $N_{V^*}$ be the corresponding operator norm on $V^*$, and let $N_{V\otimes V}$ be the operator norm $N_{V^*\to V}$, where $V\otimes V$ is naturally identified with the space of linear maps from $V^*$ to $V$ as above. We have \[N_{V\otimes V}(A)=\mopl{sup}_{N_{V^*}(\xi)=1}N_V(A\xi)=\mopl{sup}_{N_{V^*}(\xi)=1}\,\mopl{sup}_{N_{V^*}(\xi)=1}|\langle A\xi,\eta\rangle|.\]
In particular, for any $u,v\in V$ we have
\[N_{V\otimes V}(u\otimes v)=N_V(u)N_V(v).\]
The operator norm $N_{V\to V\otimes V}(\Delta)$ is a positive real number. Choosing the new norm
\[\|x\|_V:=N_{V\to V\otimes V}(\Delta)N_V(x)\]
on $V$, it is easily seen, in view of \eqref{norm-tens}, that the corresponding operator norm $\|-\|_{V\to V\otimes V}$ is obtained from $N_{V\to V\otimes V}$ by dividing it by $N_{V\to V\otimes V}(\Delta)$. In particular, it verifies $\|\Delta\|_{V\to V\otimes V}=1$. According to Proposition \ref{norm-transpose}, we also have $\|m\|_{V^*\otimes V^*\to V^*}=1$, which is equivalent to
\[
\|\xi\star\eta\|_{V^*}\le  \|\xi\|_{V^*}\|\eta\|_{V^*}
\]
 for any $\xi,\eta\in V^*$.
\end{proof}
Let us recall that all norms on a given finite-dimensional real vector space $V$ are equivalent, in the sense that, for any two norms $\|-\|_1$ and $\|-\|_2$ on $V$, there exist positive constants $k$ and $K$ such that, for any $x\in V$, the double inequality
\[k\|x\|_1\le \|x\|_2\le K\|x\|_1\]
holds. Some norms can however be more convenient than others:
\begin{definition}\label{grad-norm}
Let $V=\bigoplus_{j\ge 0}V_j$ be a graded vector space, let $\pi_j:V\to\hskip -9pt\to V_j$ be the projection to the $j$-th homogeneous component. A norm $\|-\|$ on $V$ is \textsl{graded} if $\|\pi(x)\|\le\|x\|$ for any $x\in V$ and for any finite sum $\pi$ of projections $\pi_j$.
\end{definition}
\begin{remark}\label{dual-graded-norm}
From Proposition \ref{norm-transpose}, for any finite-dimensional graded vector space $V$ endowed with a graded norm, the corresponding operator norm on $V^*$ is also graded.
\end{remark}
\section{Truncated connected graded Hopf algebras}\label{sect:truncated}
In this section we define truncated graded Hopf algebras. We follow the general Hopf algebra notations, namely $\mu,\Delta$ are used for product and  coproduct, and $u, \epsilon, S$ are used to denote unit, counit and antipode in a Hopf algebra $H$. We often write $xy$ for $\mu(x\otimes y)$.\\

Let $H=\bigoplus_{n\geq 0}H_{n}$ be a graded Hopf algebra (we include finite dimension of the homogeneous components in the definition). We can truncate it and we get $H_{\leq N} :=\bigoplus_{n=0}^{N}H_{n}$. The product and coproduct should naturally follow the truncation which can be defined with respect to the projection map $\pi:H \twoheadrightarrow H/H_{\geq N+1}$ (with self-explanatory notations). The restriction of projection map  $\pi|_{H_{\leq N}}$ is a linear isomorphism, denoted by $\phi_{N}:H_{\leq N}\rightarrow H/H_{\geq N+1}$. The product $\mu$ should therefore be defined as follows:
\begin{align*}
\mu(x\otimes y) := \phi_{N}^{-1}\big(\phi_{N}(x) \phi_{N}(y)\big),\,\text{ where }\, x,y\in  H_{\leq N}.
\end{align*}
On the other hand, the coproduct $\Delta$ takes values into
\[(H\otimes H)_{\leq N}:=\bigoplus_{p+q\le N}H_p\otimes H_q,\]
which is strictly included in $H_{\leq N}\otimes H_{\leq N}$, hence can be simply defined by restriction.
The usual compatibility diagram below does not commute,

\begin{align*}
\xymatrix{
H_{\leq N}\bigotimes H_{\leq N}\ar[r]^{\quad \mu}\ar[dd]^{\Delta \otimes \Delta} & H_{\leq N} \ar[r]^{\Delta}& H_{\leq N}\bigotimes H_{\leq N}\\
\\
(H_{\leq N})^{\otimes 4}\ar[r]^{\tau_{23}} & (H_{\leq N})^{\otimes 4}\ar[ruu]^{\mu \otimes \mu}
}
\end{align*}
but the difference $\Delta\circ\mu-(\mu\otimes\mu)\circ\tau_{23}\circ (\Delta\otimes\Delta)$ takes values into
\[(H\otimes H)_{\ge N+1}:=\bigoplus_{p+q\ge N+1}H_p\otimes H_q.\]
Following this fact we can define the notion of truncated connected Hopf algebra as follows:

\begin{definition}\label{defn:thopf}
We call $\mathcal{H}=\bigoplus_{n\geq 0}\mathcal{H}_n$ an \textsl{$N$-truncated connected graded Hopf algebra} if
\begin{itemize}
\item The $\mathcal H_n$'s are finite-dimensional, and $\mathcal{H}_n={0}$ for $n\geq N+1$.
\item The coproduct $\Delta:\mathcal{H}\rightarrow \mathcal{H}\otimes \mathcal{H}$ is a graded coassociative coproduct.
\item The $u$ and $\epsilon$ verify the axioms of a unit and a counit respectively.
\item The product $\mu:\mathcal{H} \otimes \mathcal{H}\rightarrow \mathcal{H}$ is a graded associative product.
\item The map $\Delta \circ \mu - (\mu \otimes \mu)\tau_{23}(\Delta \otimes \Delta):\mathcal{H}\otimes\mathcal{H} \rightarrow \mathcal{H}\otimes \mathcal{H}$ takes values in
\[(\mathcal{H}\otimes \mathcal{H})_{\geq N+1}=\bigoplus_{p+q\ge N+1}\mathcal H_p\otimes\mathcal H_q.\]
\end{itemize}
\end{definition}

It is easy to check that, for any connected graded Hopf algebra $H$, the subspace $\mathcal \calH:= H_{\le N}$ endowed with multiplication $\mu$ and comultiplication $\Delta$ defined above is an $N$-truncated connected graded Hopf algebra.
\begin{remark}\label{almost-comp}
As the map $\Delta \circ \mu - (\mu \otimes \mu)\tau_{23}(\Delta \otimes \Delta)$ respects the degree, its restriction to $(\mathcal H\otimes\mathcal H)_{\le N}$ is identically $0$.
 \end{remark}
\begin{remark}\label{self-dual}
The graded dual $\mathcal H^\circ$ of an $N$-truncated connected graded Hopf algebra $\mathcal{H}=\bigoplus_{n= 0}^N\mathcal{H}_n$ coincides with its full dual $\mathcal H^\star$, and is also an $N$-truncated connected graded Hopf algebra. Indeed, following the definition of a  truncated Hopf algebra we have
\begin{itemize}
\item $\mathcal H_{n}^{\circ}=0$ for any $n\geq N+1$.
\item The coproduct $\Delta_{\star} := \,^t\!\mu: \mathcal H^{\circ}\rightarrow \mathcal H^{\circ}\otimes \mathcal H^{\circ}$, where $\,^t\!\mu$ is the transpose of $\mu$, is graded coassociative.

\item $u_{\star}:=\,^t\!\epsilon: \mathbb K\rightarrow \mathcal H^{\circ}$ and $\epsilon_{\star}:=\,^t\!u: \mathcal H^{\circ} \rightarrow \mathbb K $ are unit and counit, respectively.

\item the convolution product $\star:=\, ^t\!\Delta$ is a graded associative product.
\end{itemize}
Moreover, the map
\[\Delta_{\star}\circ\star-(\star\otimes\star)\tau_{23}(\Delta_{\star}\otimes\Delta_{\star}):\mathcal H^{\circ}\otimes \mathcal H^{\circ}\rightarrow \mathcal H^{\circ}\otimes \mathcal H^{\circ}\]
takes values in $(\mathcal H^{\circ}\otimes \mathcal H^{\circ})_{\geq N+1}$. In fact, for any homogeneous $x,y \in \mathcal H$ with $|x| + |y|\leq N$ and any $\zeta, \eta\in \mathcal H^\circ$ we have
\begin{align*}
&\left\langle \Big( \Delta_{\star} \circ \star - (\star \otimes \star) \tau_{23} (\Delta_{\star} \otimes \Delta_{\star})\Big)(\zeta\otimes \eta), \,x \otimes y\right\rangle\\
=&\ \left\langle\zeta\otimes \eta,\, \Big(\Delta \circ \mu - (\mu \otimes \mu)\tau_{23}(\Delta \otimes \Delta)\Big)(x\otimes y) \right\rangle\\
=&\ 0
\end{align*}
in view of Remark \ref{almost-comp}. Hence we get
\[
\Big( \Delta_{\star} \circ \star - (\star \otimes \star) \tau_{23} (\Delta_{\star} \otimes \Delta_{\star})\Big)(\zeta\otimes \eta)\in
(\mathcal H\otimes\mathcal H)_{\le N}^\perp=(\mathcal H^\circ\otimes\mathcal H^\circ)_{\ge N+1}.
 \]
 The unit element $u_\star(1)\in \mathcal H^\circ$ coincides with the counit $\epsilon$ of $\mathcal H$, and will also be denoted by $\mathbf 1^\star$, in the same way we denote the unit element $u(1)$ of $\mathcal H$ by $\mathbf 1$.
\end{remark}

\begin{notation}\label{notation-sec-two}\rm
For later use and in coherence with Notation \ref{notation-intro}, the dimension of the homogeneous component $\mathcal H_1$ will be denoted by $d$.
\end{notation}
\section{Reminders on controls}\label{sect:controls}
This section reviews the concept of control and outlines some basic facts for later use. Let $T$ be a positive real number, and let $\Delta_T:=\{(s,t)\in[0,T]^2,\, s\le t\}$. Recall that a map $\omega:\Delta_T\to [0,+\infty[$ is a \textsl{control} if
\begin{itemize}
\item $\omega(t,t)=0$ for any $t\in[0,T]$,
\item The following superadditivity property holds for any $s\le u\le t\in[0,T]$
\begin{equation}\label{superadd}
\omega(s,u)+\omega(u,t)\le \omega(s,t).
\end{equation}
\end{itemize}

Let us present two examples.

\begin{example}\label{holder}
The map $\omega:\Delta_T\to [0,+\infty[$ given by $\omega(s,t) :=t-s$ is obviously a control, which is even additive: inequality \eqref{superadd} is an equality.
\end{example}
\begin{example}\label{pvar}
A path is a function $X$ from $[0,T]$ into a metric space$(E,d)$, not necessarily continuous. For any positive real $p$, the $p$-variation norm of $X$ is defined by
\[
\| X\|_{p,[0,T]}:=\left(\sup_{\mathcal P}\sum_{[a,b]\in \mathcal P}d(X_{a}, X_{b})^p\right)^{\frac{1}{p}},
\]
where the supremum is taken over the partitions $\mathcal P$ of the interval $[0,T]$. It is easily checked that $\omega:\Delta_T\to [0,+\infty[$ defined by
\[\omega(s,t):=\| X\|_{p,[s,t]}^p\]
is a control. Here the $p$-variation norm is defined on the interval $[s,t]$ instead of $[0,T]$.
\end{example}
\begin{proposition}\label{prop-control}\cite[Paragraph 1.1]{FZ18}
Let $\omega:[0,T]\to [0,+\infty[$ be a control. Then
\begin{enumerate}[label={\rm (\roman*)}.]
\item $\omega (s',t')\le\omega(s,t)$ for any $0\le s\le s'\le t'\le t\le T$.
\item For any $s\le t\in[0,T]$, the following limits exist
\begin{itemize}
\item $\omega(s_+,t):=\mopl{lim}_{s'\to s,\, s'>s}\omega(s',t)$.
\item $\omega(s,t_-):=\mopl{lim}_{t'\to t,\, t'<t}\omega(s',t)$.
\item $\omega(s_+,t_-):=\mopl{lim}_{\varepsilon\to 0,\, \varepsilon>0}\omega(s+\varepsilon,t-\varepsilon)$.
\item $\omega(s_-,s_+):=\mopl{lim}_{\varepsilon\to 0,\, \varepsilon>0}\omega(s-\varepsilon,s+\varepsilon)$.
\end{itemize}
\item For any $\delta>0$, the set
\[E_\delta(\omega):=\{s\in [0,T],\, \omega(s_-,s_+)\ge \delta\}\]
is finite.

\item For any $\varepsilon>0$, there exists a (finite) partition $\mathcal P$ of the interval $[0,T]$ such that, for any $[s,t]\in \mathcal P$ the majoration $\omega(s_+,t_-)<\varepsilon$ holds.
\end{enumerate}
\end{proposition}

Proposition \ref{prop-control} is actually valid for a larger class of functions called \textsl{mild controls}. This class contains any $f\circ \omega$ where $\omega$ is a control and $f:[0,+\infty[\to [0,+\infty[$ is (strictly) increasing. In particular, $\omega^\gamma$ for any $\gamma>0$ is a mild control. For a proof, see \cite[Paragraph 1.1]{FZ18}.

\begin{definition}
Let $(E',d')$ be a metric space, let $Y:[0,T]\to E'$ be a path, and let $\omega:\Delta_T\to\mathbb [0,+\infty[$ be a control. Let $\gamma>0$. We say that the path $Y$ is \textsl{$\gamma$-adapted to the control $\omega$} if, for any $(s,t)\in\Delta_T$ we have
\[d'(Y_s,Y_t)\le C\omega(s,t)^\gamma.\]
\label{defn:adapt}
\end{definition}

Let us remark that for the control $\omega(s,t)=t-s$ we recover the notion of $\gamma$-H\" older continuous path, whereas the control of Example \ref{pvar} gives back bounded $1/\gamma$-variation paths, independently of the choice of $p$:
\begin{proposition}
Let $p>0$, let $X:[0,T]\to (E,d)$ be a path with bounded $p$-variation in a metric space $(E,d)$, and let
\begin{eqnarray*}
\omega:\Delta_T&\longrightarrow [0,+\infty[\\
(s,t)&\longmapsto \|X\|_{p,[s,t]}^p
\end{eqnarray*}
be the control of Example \ref{pvar}. Let $Y:[0,T]\to (E',d')$ be a path with values in another metric space, and let $\gamma\in]0,1]$. If the path $Y$ is $\gamma$-adapted to $\omega$, it has bounded $1/\gamma$-variation.
\end{proposition}
\begin{proof}
Supposing that $Y$ is $\gamma$-adapted to $\omega$ we get
\begin{eqnarray*}
\|Y\|_{1/\gamma,[0,T]} ^{1/\gamma}&=&\mopl{sup}_{\mathcal P\smop{ partition of }[0,T]}\sum_{[a,b]\in\mathcal P}d'(Y_a,Y_b)^{1/\gamma}\\
&\le &C^{1/\gamma}\mopl{sup}_{\mathcal P\smop{ partition of }[0,T]}\sum_{[a,b]\in\mathcal P}\omega(a,b)\\
&\le&C^{1/\gamma}\mopl{sup}_{\mathcal P\smop{ partition of }[0,T]}\sum_{[a,b]\in\mathcal P}\|X\|_{p,[a,b]}^p\\
&\le&C^{1/\gamma}\mopl{sup}_{\mathcal P\smop{ partition of }[0,T]}\sum_{[a,b]\in\mathcal P}\mopl{sup}_{\mathcal Q\smop{ partition of }[a,b]}
\sum_{[x,y]\in\mathcal Q}d(X_x,Y_y)^p\\
&\le&C^{1/\gamma}\|X\|_{p,[0,T]}^p <+\infty.
\end{eqnarray*}
\end{proof}
\begin{remark}\label{controls}\rm
For any control $\omega$, and for any $\varepsilon>0$, the function $\omega^{1+\varepsilon}$ is a control as well: using the elementary inequality
\[a^{1+\varepsilon}+b^{1+\varepsilon}\le (a+b)^{1+\varepsilon}\]
for any $a,b\in \RR_{\ge 0}$, we indeed get
\[\omega(s,u)^{1+\varepsilon}+\omega(u,t)^{1+\varepsilon}\le \big(\omega(s,u)+\omega(u,t)\big)^{1+\varepsilon}\le \omega(s,t)^{1+\varepsilon}\]
for any $0\le s\le u\le t\le T$.
\end{remark}

We conclude this section with the following concept that we need.

\begin{definition}\label{rc-controls}
A control $\omega:[0,T]\to[0,+\infty[$ is \textsl{right-continuous} if for any $0\le s\le t\le T$,
\[\omega(s,t_+)=\omega(s,t).\]
\end{definition}
\noindent Note that we do not ask for right-continuity with respect to the first variable.

\section{Controlled rough paths}\label{sect:crp}
In this section, we work in the framework of truncated Hopf algebras developed in Section \ref{sect:truncated}.  From now on, we will stick to the case when $\mathcal H$ is commutative. We begin by introducing the concepts of truncated rough paths and truncated controlled rough paths. Next, we examine the composition of truncated controlled rough paths with regular functions, paying particular attention to stability. Finally, we construct the rough integral, with the aim of solving the RDE (\ref{ivp0}).
\subsection{Truncated rough paths and controlled rough paths}\label{ssect:trp}
Let $\mathcal H$ be an $N$-truncated connected graded Hopf algebra. A \textsl{truncated character} is a linear map $\phi: \mathcal H\rightarrow \mathbb R$ such that,
\begin{align*}
\phi(\mathbf 1)=1, \quad
\phi(xy)=\phi(x)\phi(y)\quad  \text{ whenever } |x|+|y|\leq N.
\end{align*}

\noindent We now introduce two key concepts that are fundamental to the present work.
\begin{remark}
From now on, the truncated Hopf algebra $\mathcal H$ is endowed with a graded norm $\|-\|$ which moreover satisfies Corollary \ref{good-norm}.
\end{remark}
\begin{definition}
Let $X:[0,T]\to \mathcal H_1^*\simeq\mathbb R^d$ be a path $\gamma$-adapted to a control $\omega$, with $\gamma\in]0,1]$ such that $N=\lfloor1/\gamma \rfloor$. In other words, $N$ is the largest integer such that $N\gamma\le 1$. An $N$-truncated rough path $\gamma$-adapted to the control $\omega$ above the path $X$ is a two-parameter family $\mathbf  X=(\mathbf  X_{st})_{s,t\in[0,T]}$ of truncated characters of the truncated Hopf algebra $\mathcal H$ verifying the following conditions:
\begin{itemize}
\item $\mathbf  X_{st}\restr{\mathcal H_1}=X_t-X_s$

\item Chen's lemma:
\[\mathbf  X_{st}=\mathbf X_{su}\star \mathbf X_{ut} \quad \text{for any $0\le s\le u\le t\le T$.}\]

\item suitable estimates:
\[
\vert\langle \mathbf  X_{st}, x \rangle\vert \leq C\|x\|_{\mathcal H}\omega(s,t)^{|x|\gamma}\quad \text{ for any homogeneous }x\in \mathcal H,
\]
where $\|-\|_{\mathcal H}$ is any graded norm (in the sense of Definition \ref{grad-norm}) on the finite-dimensional graded vector space $\calH$.
\end{itemize}
\end{definition}
\begin{definition}\label{norm-rp}
The best constant $C$ in the estimates above is the norm of the truncated rough path $\mathbf X$, namely
\[\|\mathbf X\|_{\omega,\gamma}:=\mopl{sup}_{j=0,\ldots N}\mopl{sup}_{(s,t)\in \Delta_T}\|\pi_j^*(\mathbf X_{st})\|_{\mathcal H^*}\omega(s,t)^{-j\gamma},\]
where $\pi_j^*:\mathcal H^*\to\hskip -9pt\to \mathcal H_j^*$ is the projection onto the $j$-th homogeneous component.
\end{definition}
Now let $\mathcal H$ be a connected graded commutative $N$-truncated Hopf algebra.
Let $\omega$ be a control, let $\gamma>0$, and let $\mathbf  X$ be an $N$-truncated rough path associated to $\mathcal H$, that we suppose $\gamma$-adapted to the control $\omega$.
\begin{definition}\label{def:crp}
An $\mathbf  X$-controlled rough path is a map
\[\mathbf  Z=(\mathbf  Z^1,\ldots,\mathbf  Z^e):[0,T]\longrightarrow \mathcal H^e\]
such that
\begin{equation}\label{estimate-crp}
\| \langle h,R_{\mathbf X}(\mathbf Z)_{st}\rangle\|_{\mathbb R^e} \le C\|h\|_{\mathcal H^*}\omega(s,t)^{(N-\vert h\vert)\gamma}\quad \text{ for any } s\le t\in [0,T]\text{ and any homogeneous } h\in \mathcal H^\star,
\end{equation}
where the remainder $R_{\mathbf X}(\mathbf Z): \Delta_T\to \calH^e$ is given by
\begin{equation}\label{crp}
\langle h,R_{\mathbf X}(\mathbf Z)_{st}\rangle:= \langle h,\mathbf  Z_t\rangle - \langle \mathbf  X_{st}\star h,\mathbf  Z_s\rangle\quad \text{ for any }  h\in \mathcal H^\star.
\end{equation}
\end{definition}

Controlled rough paths are organized into a Banach space: our approach is close to \cite[Paragraph 4.3]{FZ18}. Let $N$ be a positive integer. For any $N$-truncated graded connected Hopf algebra $\mathcal H$ and for any $j\in\{0,\ldots, N\}$, let $\pi_j:\mathcal H^e\to\hskip -9pt\to \mathcal H_j^e$ be the componentwise projection of the $j$-th homogeneous part. In view of \eqref{estimate-crp}, for any map $R:\Delta_T\to \mathcal H^e$ we set
\[\|R\|_{\omega,\gamma}:=\mopl{sup}_{j\in\{0,\ldots,N\}}\, \mopl{sup}_{(s,t)\in\Delta_T}\|\pi_j (R_{st})\|_{\mathcal H^e}\, \omega(s,t)^{-(N-j)\gamma}.\]

Let $\gamma\in ]0,1]$. For any truncated rough path $\mathbf X$, supposed to be $\gamma$-adapted to a control $\omega$, we denote by $\mathcal V_{\mathbf X}^{[0,T]}$ the space of controlled rough paths in the sense of Definition \ref{def:crp}. The controlled rough path norm
\begin{equation}\label{crp-norm}
\|\mathbf Z\|:=\|\mathbf Z_0\|_{\mathcal H^e}+\|R_{\mathbf X}(\mathbf Z)\|_{\omega,\gamma}
\end{equation}
 is finite for any $\mathbf Z\in \mathcal V_{\mathbf X}^{[0,T]}$, and turns $\mathcal V_{\mathbf X}^{[0,T]}$ into a Banach space. Let us emphasize at this stage that $\|\mathbf Z\|$ is a nondecreasing function of the upper end $T$ of the interval $[0,T]$ on which the controlled rough path is considered. In other words, for any $T'\in]0,T]$, the restriction of controlled rough paths to $[0,T']$ is a weak contraction from $\mathcal V_{\mathbf X}^{[0,T]}$ onto $\mathcal V_{\mathbf X}^{[0,T']}$ in the sense that, for any controlled rough path $\mathbf Z\in \mathcal V_{\mathbf X}^{[0,T]}$, we have
\begin{equation*}
\|\mathbf Z\restr{[0,T']}\|\le \|\mathbf Z\|.
\end{equation*}
\begin{remark}\label{rmk-deltaZ}
The path $Z:[0,T]\to \mathbb R^e$ defined by
$$Z_t:=\langle \mathbf 1^\star,\,\mathbf  Z_t\rangle=\left(\langle \mathbf 1^\star,\,\mathbf  Z^1_t\rangle,\ldots, \langle \mathbf 1^\star,\,\mathbf  Z^e_t\rangle\right)$$
verifies the following estimate:
\begin{equation}\label{z control}
\|Z_{t} - Z_{s}\|_{\mathbb R^e} \leq  \mop{max}\big(1,\omega(0,T)\big)^{(N-1)\gamma}\left(N\|\mathbf X\|_{\omega,\gamma}\,\|\mathbf Z\|_{\smop{sup}}+\|\mathbf Z\|\right)\,\omega(s,t)^{\gamma}.
\end{equation}
where $\|\mathbf Z\|_{\smop{sup}}:=\mopl{sup}_{t\in[0,T]}\|\mathbf Z_t\|_{\mathcal H^e}$. Indeed,
\begin{eqnarray*}
\|Z_t-Z_s\|_{\mathbb R^e}&=&\|\langle \mathbf 1^*,\, \mathbf Z_t-\mathbf Z_s\rangle\|\\
&\le & \|\langle\mathbf X_{st}-\mathbf 1^*,\, \mathbf Z_s\rangle\|+\|\langle\mathbf 1^*,\, R_{\mathbf X}(\mathbf Z)_{st}\rangle\|\\
&\le &\|\langle \mathbf X_{st},\, \mathbf Z_s-Z_s\mathbf 1\rangle\|+\|\langle\mathbf 1^*,\, R_{\mathbf X}(\mathbf Z)_{st}\rangle\|\\
&\le &\sum_{j=1}^N\|\langle \mathbf X_{st},\, \pi_j(\mathbf Z_s)\rangle\|+\|\langle\mathbf 1^*,\, R_{\mathbf X}(\mathbf Z)_{st}\rangle\|\\
&\le &\sum_{j=1}^N \|\mathbf X\|_{\omega,\gamma}\|\mathbf Z_s\|\omega(s,t)^{j\gamma}+\|R_{\mathbf X}(\mathbf Z)\|_{\omega,\gamma}\omega(s,t)^{N\gamma}\\
&\le & \mop{max}\big(1,\omega(0,T)\big)^{(N-1)\gamma}\left(N\|\mathbf X\|_{\omega,\gamma}\,\|\mathbf Z\|_{\smop{sup}}+\|\mathbf Z\|\right)\,\omega(s,t)^{\gamma}.
\end{eqnarray*}
\end{remark}

\begin{proposition}\label{fz}(see (31) in \cite{FZ18})
For any controlled rough path $\mathbf Z$ driven by the rough path $\mathbf X$, its supremum norm (see Remark \ref{rmk-deltaZ}) verifies
\begin{equation}\label{sup-crp}
\|\mathbf Z\|_{\smop{sup}}\le (N+1)K(T)^N\left(1+\|\mathbf X\|_{\omega,\gamma}\right)\|\mathbf Z\|,
\end{equation}
with $K(T):=\mop{max}\big(1,\omega(0,T)\big)^{\gamma}$.
\end{proposition}
\begin{proof}
From the definition of the remainder, we have for any $(s,t)\in \Delta_T$:
\[\mathbf Z_t=R(\mathbf Z)_{st}+T_{st}(\mathbf Z_s),\]
where $T_{st}$ is (the coordinate-wise action on $\mathcal H^e$ of) the transpose of the convolution operator $\mathbf X_{st}\star -$. From Corollary \ref{good-norm} we have
\[\|T_{st}\|_{\mathcal H^e\to\mathcal H^e}\le \|\mathbf X_{st}\|_{\mathcal H^*}.\]
We therefore have
\begin{eqnarray*}
\|\mathbf Z\|_{\smop{sup}}&=&\mopl{sup}_{t\in[0,T]}\|\mathbf Z_t\|\\
&\le &\mopl{sup}_{t\in[0,T]}\|R(\mathbf Z)_{0t}\|+\big(\mopl{sup}_{t\in[0,T]}\|\mathbf X_{0t}\|\big)\,\|\mathbf Z_0\|\\
&\le &\mopl{sup}_{(s,t)\in\Delta_T}\|R(\mathbf Z)_{st}\|+\big(\mopl{sup}_{(s,t)\in\Delta_T}\|\mathbf X_{st}\|\big)\,\|\mathbf Z_0\|\\
&\le&\sum_{j=0}^N\mopl{sup}_{(s,t)\in\Delta_T}\|\pi_j(R(\mathbf Z)_{st)}\|+\sum_{j=0}^N\big(\mopl{sup}_{(s,t)\in\Delta_T}\|\pi_j^*(\mathbf X_{st})\|\big)\,\|\mathbf Z_0\|\\
&\le &\sum_{j=0}^N\|R(\mathbf Z)\|_{\omega,\gamma}\mop{max}\big(1,\omega(0,T)\big)^{(N-j)\gamma}+\|\mathbf Z_0\|\sum_{j=0}^N\|\mathbf X\|_{\omega,\gamma}\mop{max}\big(1,\omega(0,T)\big)^{j\gamma}\\
&\le &(N+1)K(T)^N\left(\|\mathbf Z\|+\|\mathbf Z_0\|.\|\mathbf X\|_{\omega,\gamma}\right) \hspace{2cm} \hbox{ (from \eqref{crp-norm})}\\
&\le & (N+1)K(T)^N\left(1+\|\mathbf X\|_{\omega,\gamma}\right)\|\mathbf Z\|.
\end{eqnarray*}
\end{proof}

\begin{corollary}
Under the setting of Proposition~\ref{fz}, the following estimate holds
\begin{equation}\label{fz-bis}
\|Z_t-Z_s\|_{\mathbb R^e}\le C\left(T,\|\mathbf X\|_{\omega,\gamma}\right)\omega(s,t)^\gamma\|\mathbf Z\|,
\end{equation}
with
\begin{equation}\label{degree-two}
C(T,x):=K(T)^{N-1}\left(N(N+1)x(x+1)K(T)^N+1\right).
\end{equation}
\end{corollary}
\begin{proof}
It is a directly checked by plugging \eqref{sup-crp} into \eqref{z control}.
\end{proof}
\subsection{Image by a smooth map}\label{ssect:phi}
Let $\varphi$ be in the Fr\'echet space $C^N(\mathbb R^e,\mathbb R^e)$ endowed with the semi-norms
\begin{equation}\label{seminorms}
\|\varphi\|^N_K:=\sum_{j=0}^N\left\|D^j\varphi(z)\right\|_{S^j(\mathbb R^e)\to\mathbb R^e,\, z\in K}.
\end{equation}
Here, $D^j\varphi:\mathbb R^e\to\mathcal L\big(S^j(\mathbb R^e)\big),\mathbb R^e)$ denotes the $j$-th differential of $\varphi$, and $K$ is a compact subset of $\mathbb R^e$.
An $\mathcal H^e$-valued path $\varphi(\mathbf  Z)$ above the $\mathbb R^e$-valued path $\varphi\circ Z$ is given through its evaluation against any $h\in \mathcal H^\star$ by the following formulae:
\begin{eqnarray*}
\langle \mathbf 1^\star, \varphi(\mathbf  Z)_s\rangle&:=&\varphi(Z_s),\\
\langle h, \varphi(\mathbf  Z)_s\rangle&:=&\sum_{n=1}^{N}\frac{D^n\varphi(Z_s)}{n!} \left\langle \Delta_\star^{\smop{red}(n-1)}(h),\,\mathbf  Z_s^{\otimes n}\right\rangle \hbox{ for any }h\in\mop{Ker}\epsilon_\star.
\end{eqnarray*}
Recall that the reduced coproduct defined by $$\Delta_\star^{\smop{red}}(h):=\Delta_\star(h)-\mathbf 1^\star\otimes h-h\otimes \mathbf 1^\star$$ is also coassociative, and that $\mathbf Z$ is a controlled rough path above $Z$. The map
$$\Delta_\star ^{\smop{red}(n-1)}:\mathcal H^\star\to(\mathcal H^\star)^{\otimes n}$$
stands for the iterated reduced coproduct, recursively defined by
\[\Delta_\star ^{\smop{red}(0)}(h):=h,\hskip 12mm \Delta_\star ^{\smop{red}(n-1)}(h):=(\Delta_\star^{\smop{red}}\otimes \mop{Id}_{\mathcal H}^{\otimes (n-2)})\circ\Delta_\star ^{\smop{red}(n-2)}(h)\quad \text{for any $h\in H^{\star}$.}\]
It vanishes identically whenever $n\ge N$.

\begin{lemma}\label{red-grp}
If $h$ is group-like for $\Delta_\star$, then $h-\mathbf 1$ is group-like for $\Delta_\star^{\smop{red}}$.
\ignore{\color{blue}{
\begin{align*}
\Delta_\star^{\smop{red}}(h-\mathbf 1)=& \Delta_\star(h-\mathbf 1) - (h-\mathbf 1)\otimes \mathbf 1 - \mathbf 1\otimes (h-\mathbf 1)\\
=& \Delta_\star(h) - \Delta_\star(\mathbf 1)  - (h-\mathbf 1)\otimes \mathbf 1 - \mathbf 1\otimes (h-\mathbf 1)\\
=& h\otimes h - \mathbf 1\otimes \mathbf 1  - (h-\mathbf 1)\otimes \mathbf 1 - \mathbf 1\otimes (h-\mathbf 1)\\
=& (h - \mathbf 1)\otimes (h - \mathbf 1).
\end{align*}
}
}
\end{lemma}

\begin{proof}
Checking $\Delta_\star^{\smop{red}}(h-\mathbf 1)=(h-\mathbf 1)\otimes(h-\mathbf 1)$ is obvious and left to the reader.
\end{proof}

As anticipated, applying a smooth map to a controlled rough path yields another controlled rough path.

\begin{proposition}
The $\mathcal H^e$-valued path $\varphi(\mathbf  Z)$ is a controlled rough path above $\varphi(Z)$, i.e. its remainder $R_{\mathbf X}\big(\varphi(\mathbf Z)\big)$ defined by \eqref{crp}, namely
\begin{equation*}
\left\langle h,R_{\mathbf X}\big(\varphi(\mathbf Z)\big)_{st}\right\rangle:=\langle h,\,\varphi(\mathbf  Z)_t\rangle-\langle \mathbf  X_{st}\star h,\,\varphi(\mathbf  Z)_s\rangle
\end{equation*}
for any $h\in\mathcal H^\star$, verifies \eqref{estimate-crp}.
\end{proposition}

\begin{proof}
Let us start with the case $h=\mathbf 1^*$. The Taylor formula for
$$\psi:u\mapsto \frac{t-u}{t-s}\varphi(Z_s)+\frac{u-s}{t-s}\varphi(Z_t)$$ at $u=s$ yields
\begin{eqnarray*}
\langle \mathbf 1^*,\, \varphi(\mathbf  Z)_t\rangle &=&\varphi(Z_t)\\
&=&\sum_{n=0}^{N-1} \frac{D^n\varphi(Z_s)}{n!}\big((Z_t-Z_s)^n\big)+R_{st}\\
&=&\varphi(Z_s)+\sum_{n=1}^{N-1}\frac{D^n\varphi(Z_s)}{n!}\big((Z_t-Z_s)^n\big)+R_{st}
\end{eqnarray*}
with $$\|R_{st}\|\le C\|Z_t-Z_s\|^{N}\le C\omega(s,t)^{N\gamma},$$
thanks to Hypothesis \eqref{z control} on the path $Z$. Here, $(Z_t-Z_s)^n$ is to be understood as an element in the $n$-th symmetric power $S^n(\mathbb R^e)$.
Now,
\allowdisplaybreaks{
\begin{eqnarray*}
\langle \mathbf  X_{st},\, \varphi(\mathbf  Z)_s\rangle &=&\varphi(Z_s)+\langle \mathbf  X_{st}-\mathbf 1,\, \varphi(\mathbf  Z)_s\rangle\\
&=&\varphi(Z_s)+\sum_{n=1}^{N}\frac{D^n\varphi(Z_s)}{n!}\left\langle \Delta_\star^{\smop{red}(n-1)}(\mathbf  X_{st}-\mathbf 1),\,\mathbf  Z_s^{\otimes n}\right\rangle\\
&=&\varphi(Z_s)+\sum_{n=1}^{N}\frac{D^n\varphi(Z_s)}{n!}\left\langle (\mathbf  X_{st}-\mathbf 1)^{\otimes n},\,\mathbf  Z_s^{\otimes n}\right\rangle \quad\quad \hbox{ (from Lemma \ref{red-grp})}\\
&=&\varphi(Z_s)+\sum_{n=1}^{N-1}\frac{D^n\varphi(Z_s)}{n!}\big((Z_t-Z_s)^n\big)+\frac{D^N\varphi(Z_s)}{N!}\big((Z_t-Z_s)^N\big)+O\left(\omega(s,t)^{N\gamma}\right),
\end{eqnarray*}
}
which ends up the case $h=\mathbf 1^*$. The case $h\in\mop{Ker}\epsilon^*$ is handled as follows (with $\overline x:=x-\epsilon(x)\mathbf 1$ for any $x\in\mathcal H$, extended component-wisely to $x\in\mathcal H^e$):
\allowdisplaybreaks{
\begin{eqnarray*}
&&\langle \mathbf  X_{st}\star h,\, \varphi(\mathbf  Z)_s\rangle \\
&=&\sum_{n=1}^{N}\frac{D^n\varphi(Z_s)}{n!}\left\langle\Delta_\star^{\smop{red}(n-1)}(\mathbf  X_{st}\star h),\,\mathbf  Z_s^{\otimes n}\right\rangle\\
&=&\sum_{n=1}^{N}\frac{D^n\varphi(Z_s)}{n!}\left\langle\Delta_\star^{(n-1)}(\mathbf  X_{st}\star h),\,\overline{\mathbf  Z_s}^{\,\otimes n}\right\rangle\\
&=&\sum_{n=1}^{N}\frac{D^n\varphi(Z_s)}{n!}\left\langle\mathbf  X_{st}^{\otimes n}\star\Delta_\star^{(n-1)}(h),\,\overline{\mathbf  Z_s}^{\,\otimes n}\right\rangle\\
&=&\sum_{n=1}^{N}\frac{D^n\varphi(Z_s)}{n!}\sum_{(h)}\left\langle \mathbf  X_{st}\star h_1,\, \overline{\mathbf  Z_s}\right\rangle\cdots\left\langle \mathbf  X_{st}\star h_n,\, \overline{\mathbf  Z_s}\right\rangle\\
&=&\hskip -3mm\sum_{n=1}^{N}\frac{D^n\varphi(Z_s)}{n!}\sum_{(h)}\Big(\langle h_1,\, \overline{\mathbf  Z_t}\rangle+(Z_t-Z_s)\langle h_1,\, \mathbf 1\rangle\Big)\cdots \Big(\langle h_n,\, \overline{\mathbf  Z_t}\rangle+(Z_t-Z_s)\langle h_n,\, \mathbf 1\rangle\Big)-R,
\end{eqnarray*}
}
where, from from \eqref{crp}, we have $R=O\left(\omega(s,t)^{(N-|h|)\gamma}\right)$. Proceeding with the computation, we have
\begin{eqnarray*}
\langle \mathbf  X_{st}\star h,\, \varphi(\mathbf  Z)_s\rangle &=&\sum_{n=1}^{N}\frac{D^n\varphi(Z_s)}{n!}\sum_{p=1}^n{n\choose p}\left\langle\Delta_\star^{\smop{red}(p-1)}(h),\, \mathbf  Z_t^{\otimes p}\right\rangle (Z_t-Z_s)^{n-p}-R\\
&=&\sum_{p=1}^{N}\sum_{q=0}^{N-p}\frac{D^{p+q}\varphi(Z_s)}{p!q!}\left\langle\Delta_\star^{\smop{red}(p-1)}(h),\, \mathbf  Z_t^{\otimes p}\right\rangle (Z_t-Z_s)^{q}-R.
\end{eqnarray*}
The inner sum is the Taylor expansion at $s$ of
\[t\mapsto \frac{D^p\varphi(Z_t)}{p!}\left\langle\Delta_\star^{\smop{red}(p-1)}(h),\, \mathbf  Z_t^{\otimes p}\right\rangle\]
at order $N-p$. We therefore get
\begin{eqnarray*}
\langle \mathbf  X_{st}\star h,\, \varphi(\mathbf  Z)_s\rangle&=&\sum_{p=1}^{N}\frac{D^p\varphi(Z_t)}{p!}\left\langle\Delta_\star^{\smop{red}(p-1)}(h),\, \mathbf  Z_t^{\otimes p}\right\rangle-R\\
&=&\langle h,\,\varphi(\mathbf  Z)_t\rangle -R.
\end{eqnarray*}
This completes the proof. From \eqref{crp}, the term $R$ in the computation is the remainder $R_{\mathbf X}\big(\varphi(\mathbf Z)\big)$ of the controlled path $\varphi\big(\mathbf Z)$.
\end{proof}

We are now ready to examine the continuity of the map $\mathbf Z\mapsto \varphi (\mathbf Z)$ in the Banach space $\mathcal V^{\mathbf X}_{[0,T]}$ of $\mathbf X$-driven controlled rough paths.

\begin{proposition}\label{main-est-varphi}
Let $\mathbf Z,\,\mathbf Z'$ be two controlled rough paths from $[0,T]$ to $\mathcal H^e$, driven by a rough path $\mathbf X$\ignore{such that $\mathbf Z_0=\mathbf Z'_0$, driven by the same rough path $\mathbf X$}, and let $\varphi\in C^N(\mathbb R^e,\,\mathbb R^e)$. There exists a constant $C_{\varphi}(T)$ such that
\[\|\varphi(\mathbf Z')-\varphi(\mathbf Z)\|\le C_{\varphi}(T)\|\mathbf Z'-\mathbf Z\|,\]
which can be chosen non-decreasing with respect to $T$. It depends on $\|\mathbf X\|_{\omega,\gamma},\,\|\mathbf Z\|, \|\mathbf Z'\|$ and $\|\varphi\|^N_{K(Z,Z')}$, where $K(Z,Z')$ is the compact subset $Z([0,T])\cup Z'([0,T])$ of $\mathbb R^e$.
\end{proposition}
\begin{proof}
According to the computation in last proposition we get
\begin{eqnarray*}
\langle \mathbf X_{st} \star h, \varphi (\mathbf Z)_s\rangle &=& \sum_{n=0}^{N} \frac{D^{n}\varphi(Z_s)}{n!}\left\langle \Delta_{\star}^{(n-1)}(\mathbf X_{st}\star h), \overline {\mathbf Z_s}^{\,\otimes n}\right\rangle\\
&=& \sum_{n=0}^{N}\frac{D^{n}\varphi(Z_{s})}{n!}\left\langle \mathbf X_{st}^{\otimes n}\star \Delta_{\star}^{(n-1)}(h), \overline{\mathbf Z_s}^{\, \otimes n}\right\rangle\\
&=&\sum_{n=1}^{N}\frac{D^n\varphi(Z_s)}{n!}\sum_{(h)}\left\langle \mathbf  X_{st}\star h_1,\, \overline{\mathbf  Z_s}\right\rangle\cdots\left\langle \mathbf  X_{st}\star h_n,\, \overline{\mathbf  Z_s}\right\rangle\\
&=& \sum_{n=0}^{N}\frac{D^{n}\varphi(Z_s)}{n!}\sum_{(h)}\prod_{j=1}^{n}\langle \mathbf X_{st}\star h_j, \mathbf Z_s - Z_s\mathbf 1\rangle.
\end{eqnarray*}
From \eqref{crp} we therefore get
\allowdisplaybreaks{
\begin{eqnarray}
&&\langle \mathbf X_{st} \star h, \varphi (\mathbf Z)_s\rangle\nonumber\\
&=& \sum_{n=0}^{N}\frac{D^{n}\varphi(Z_s)}{n!}\sum_{(h)}\prod_{j=1}^{n}\left(\langle h_j, \mathbf Z_t\rangle - \langle h_j, R_{\mathbf X}(\mathbf Z)_{st} \rangle - \langle h_j, Z_s\mathbf 1\rangle\right)\nonumber\\
&=& \sum_{n=0}^{N}\frac{D^{n}\varphi(Z_s)}{n!}\sum_{(h)}\prod_{j=1}^{n}\left(\langle h_j, \overline{\mathbf Z_t}\rangle + (Z_t - Z_s)\langle h_j, \mathbf 1\rangle - \langle h_j, R_{\mathbf X}(\mathbf Z)_{st} \rangle\right)\nonumber\\
&=& \sum_{n=0}^{N}\frac{D^{n}\varphi(Z_s)}{n!}\sum_{(h)}\sum_{I\subseteq \{1,\ldots, n\}}(-1)^{n-|I|}\prod_{j\in I}\left(\langle h_j, \overline{\mathbf Z_t}\rangle + (Z_t - Z_s)\langle h_j, \mathbf 1 \rangle\right)\prod_{j\notin I}\langle h_j, R_{\mathbf X}(\mathbf Z)_{st} \rangle\nonumber\\
&=& \sum_{n=0}^{N}\frac{D^{n}\varphi(Z_s)}{n!}\sum_{(h)}\prod_{j=1}^{n}\left(\langle h_j, \overline{\mathbf Z_t}\rangle + (Z_t - Z_s)\langle h_j, \mathbf 1 \rangle\right)\nonumber\\
& &  + \sum_{n=0}^{N}\frac{D^{n}\varphi(Z_s)}{n!}\sum_{(h)}\sum_{I\subsetneq\{1,\ldots,n\}}(-1)^{n-|I|}\prod_{j\in I}\left(\langle h_j, \overline{\mathbf Z_t}\rangle + (Z_t - Z_s)\langle h_j, \mathbf 1\rangle\right)\prod_{j\notin I}\langle h_j, R_{\mathbf X}(\mathbf Z)_{st}\rangle\nonumber\\
&=& \sum_{n=0}^{N}\frac{D^{n}\varphi(Z_s)}{n!}\sum_{(h)}\prod_{j=1}^{n}\left(\langle h_j, \overline{\mathbf Z_t}\rangle + (Z_t - Z_s)\langle h_j, \mathbf 1 \rangle\right)\label{term-one}\\
& & + \sum_{n=0}^{N}\frac{D^{n}\varphi(Z_s)}{n!}\sum_{(h)}\sum_{I\subsetneq\{1,\ldots,n\}}(-1)^{n-|I|}\prod_{j\in I}\left(\langle h_j, \mathbf Z_t - Z_s\mathbf 1\rangle\right)\prod_{j\notin I}\langle h_j, R_{\mathbf X}(\mathbf Z)_{st}\rangle\label{term-two}\\
&=& \left\langle h, \varphi(\mathbf Z)_t\right\rangle + \left\langle h, R_{\mathbf X}\big(\varphi(\mathbf Z)\big)_{st}\right\rangle,\nonumber
\end{eqnarray}
}
where $\left\langle h, \varphi(\mathbf Z)_t\right\rangle$ coincides with \eqref{term-one}, and $ \left\langle h, R_{\mathbf X}\big(\varphi(\mathbf Z)\big)_{st}\right\rangle$ with \eqref{term-two}. A similar computation of course holds for $\mathbf Z'$. In order to explain the details we define the elements
$\widetilde{\mathbf z}_j^{I}$ and
$\widetilde{\mathbf z}_j^{\, 'I}$ as follows:
\[
\widetilde{\mathbf z}_j^I :=
\setlength{\abovedisplayskip}{5pt}
\setlength{\belowdisplayskip}{5pt}
\begin{cases}
\mathbf Z_t - Z_s\mathbf1, \ &\text{ if } j\in I\\
R_{\mathbf X}(\mathbf Z)_{st}, \ &\text{ if } j \notin I,
\end{cases}
\hskip 12mm
\widetilde{\mathbf z}_j^{\, 'I} :=
\setlength{\abovedisplayskip}{5pt}
\setlength{\belowdisplayskip}{5pt}
\begin{cases}
\mathbf Z'_t - Z'_s\mathbf1, \ &\text{ if } j\in I\\
R_{\mathbf X}(\mathbf Z')_{st}, \ &\text{ if } j \notin I.
\end{cases}
\]
The computation above yields
\begin{eqnarray*}
\| \langle h, R_{\mathbf X}\big(\varphi(\mathbf Z')\big)_{st} - R_{\mathbf X}\big( \varphi(\mathbf Z)\big)_{st}\rangle\|&\le& \bigg\| \sum_{n=0}^N \frac{D^n\varphi(Z_s)}{n!} \sum_{(h)}\sum_{I \subsetneq \{1,\ldots,n\}}(-1)^{n-\vert I \vert}\prod_{j=1,\ldots,n}\langle h_j,  \widetilde{\mathbf z}_j^{I}\rangle\\
&&- \sum_{n=0}^N \frac{D^n\varphi(Z'_s)}{n!} \sum_{(h)}\sum_{I \subsetneq \{1,\ldots,n\}}(-1)^{n-\vert I \vert}\prod_{j=1,\ldots,n}\langle h_j,  \widetilde{\mathbf z}_j^{\, 'I}\rangle\bigg\|.
\end{eqnarray*}
Using the equality
\begin{eqnarray*}
a_0'\cdots a_n' - a_0\cdots a_n & = & \sum_{j=0}^{n} a_0'\cdots a_{j-1}'(a_j' - a_j)a_{j+1}\cdots a_n
\end{eqnarray*}
valid in any associative algebra, we therefore get
\begin{eqnarray*}
\left\| \langle h, R_{\mathbf X}\big(\varphi(\mathbf Z')\big)_{st} - R_{\mathbf X}\big( \varphi(\mathbf Z)\big)_{st}\rangle\right\|\hskip -25mm &&\\
&\le& \left\| \sum_{n=0}^N \frac{D^n\varphi(Z_s)}{n!} \sum_{(h)}\sum_{I \subsetneq \{1,\ldots,n\}}(-1)^{n-\vert I \vert}\sum_{j=1}^n\langle h_1, \widetilde{\mathbf z}_1^{\, 'I}\rangle \cdots \langle h_j, \widetilde{\mathbf z}_j^{\, 'I} - \widetilde{\mathbf z}_j^{I}\rangle \cdots \langle h_n, \widetilde{\mathbf z}_n^{I}\rangle \right\|\\
&&+ \left\|\frac{D^n\varphi(Z'_s)-D^n\varphi(Z_s)}{n!} \sum_{(h)}\sum_{I \subsetneq \{1,\ldots,n\}}(-1)^{n-\vert I \vert}\langle h_1, \widetilde{\mathbf z}_1^{I}\rangle \cdots \langle h_n, \widetilde{\mathbf z}_n^{I}\rangle \right\|\\
& \leq & \sum_{n=0}^N \frac{\| D^n\varphi(Z_s)\|}{n!}\sum_{(h)}\sum_{j=1}^n\sum_{I\subsetneq \{1,\ldots, n\}}\| h_1\|\cdots\| h_n\| \ \|\widetilde{\mathbf z}_1^{\, 'I}\| \cdots \| \widetilde{\mathbf z}_j^{\, 'I} - \widetilde{\mathbf z}_j^{I}\| \cdots \| \widetilde{\mathbf z}_n^{I}\|\\
&&+\sum_{n=0}^N \frac{\| D^n\varphi(Z_s)-D^n\varphi(Z'_s)\|}{n!}\sum_{(h)}\sum_{I\subsetneq \{1,\ldots, n\}}\| h_1\|\cdots\| h_n\| \ \|\widetilde{\mathbf z}_1^{I}\| \cdots \| \widetilde{\mathbf z}_n^{I}\|\\
& \le & \|h\|\,\sum_{n=0}^N \frac{\| D^n\varphi(Z_s)\|}{n!}\,\|\Delta_{\star}^{(n-1)}\|\sum_{j=1}^n\sum_{I\subsetneq \{1,\ldots, n\}} \ \| \widetilde{\mathbf z}_1^{\, 'I}\| \cdots \| \widetilde{\mathbf z}_j^{\, 'I} - \widetilde{\mathbf z}_j^{I}\| \cdots \| \widetilde{\mathbf z}_n^{I}\|\\
&&+\|h\|\,\sum_{n=0}^N \frac{\| D^n\varphi(Z_s)-D^n\varphi(Z'_s)\|}{n!}\,\|\Delta_{\star}^{(n-1)}\|\sum_{I\subsetneq \{1,\ldots, n\}} \ \| \widetilde{\mathbf z}_1^{I}\| \cdots \| \widetilde{\mathbf z}_n^{I}\|.
\end{eqnarray*}\label{est-one}

\noindent We easily check the estimate
\begin{align*}
\|\mathbf Z_t-Z_s\mathbf 1\|&\le \|\mathbf Z_t-Z_s\mathbf 1\|+\|(Z_t-Z_s)\mathbf 1\|\le \|\mathbf Z_t\|+\|Z_t-Z_s\|_{\mathbb R^e},
\end{align*}
hence
\[\|\mathbf Z_t-Z_s\mathbf 1\|\le\bigg((N+1)K(T)^N\left(1+\|\mathbf X\|_{\omega,\gamma}\right)+C\left(T,\|\mathbf X\|_{\omega,\gamma}\right)\omega(0,T)^\gamma\bigg) \|\mathbf Z\|\]
according to \eqref{sup-crp} and \eqref{fz-bis}.\\

\noindent Supposing $h\in \mathcal H_r$, and using the component-wise projection map
\[\pi_r:\mathcal H^e\to\hskip -9pt \to \mathcal H_r^e,\]
we therefore have the following estimate of the remainder:
\begin{eqnarray*}
&& \left\|\pi_r \Big(R_{\mathbf X}\big(\varphi(\mathbf Z^{'})\big)_{st} - R_{\mathbf X}\big( \varphi(\mathbf Z)\big)_{st}\Big)\right\| \nonumber\\
 &\le& \sum_{n=0}^{N} \frac{\|D^n\varphi(Z_s)\|}{n!} \|\Delta_\star^{(n-1)}\|\Bigg\{\sum_{I \subsetneq \{1,\ldots, n\}} \vert I\vert\, \bigg(\mop{max}\left(\big \| \pi_r ( \mathbf Z_t)\big\|,\big\|\pi_r(\mathbf Z_t^{'})\big\|\right)+ C\left(T,\|\mathbf X\|_{\omega,\gamma}\right)\omega(s,t)^{\gamma}\bigg)^{\vert I\vert-1}\nonumber\\
&& \hskip 25mm \times\left\| \pi_r\left(\mathbf Z_t^{'} - \mathbf Z_t + (Z'_s-Z_s)\mathbf 1\right)\right\|\, \mop{max}\Big(\left\|\pi_rR_{\mathbf X}(\mathbf Z)_{st}\right\|,\,\left\|\pi_rR_{\mathbf X}(\mathbf Z')_{st}\right\|\Big)^{n - \vert I\vert}\nonumber\\
& & \hskip 2mm+\sum_{I \subsetneq \{1,\ldots, n\}} (n-\vert I\vert)\, \bigg(\mop{max}\left(\big \| \pi_r ( \mathbf Z_t)\big\|,\big\|\pi_r(\mathbf Z_t^{'})\big\|\right) + C\left(T,\|\mathbf X\|_{\omega,\gamma}\right)\omega(s,t)^{\gamma}\bigg)^{\vert I\vert}\nonumber\\
&&  \hskip 25mm \times \left\| \pi_r\big(R_{\mathbf X}(\mathbf Z^{'})_{st} - R_{\mathbf X}(\mathbf Z)_{st}\big)\right\| \mop{max}\Big(\left\|\pi_rR_{\mathbf X}(\mathbf Z)_{st}\right\|,\,\left\|\pi_rR_{\mathbf X}(\mathbf Z')_{st}\right\|\Big)^{n-\vert I\vert - 1}\Bigg\}\\
&& \hskip 2mm +\sum_{n=0}^{N} \frac{\|D^n\varphi(Z'_s)-D^n\varphi(Z_s)\|}{n!}\sum_{I\subsetneq\{1,\ldots,n\}}|I| \bigg(\mop{max}\left(\big \| \pi_r ( \mathbf Z_t)\big\|,\big\|\pi_r(\mathbf Z_t^{'})\big\|\right)+ C\left(T,\|\mathbf X\|_{\omega,\gamma}\right)\omega(s,t)^{\gamma}\bigg)^{|I|}\\
&&\hskip 25 mm \times \mop{max}\Big(\left\|\pi_rR_{\mathbf X}(\mathbf Z)_{st}\right\|,\,\left\|\pi_rR_{\mathbf X}(\mathbf Z')_{st}\right\|\Big)^{n-|I|}.
\end{eqnarray*}
\noindent After dividing by $\omega(s,t)^{(N-r)\gamma}$ on both sides, we have
\begin{eqnarray}\label{est-two}
&& \left\| \pi_r \Big(R_{\mathbf X}\big(\varphi(\mathbf Z')\big)_{st}-R_{\mathbf X}\big(\varphi(\mathbf Z)\big)_{st}\Big)\right\|\omega(s,t)^{-(N-r)\gamma}\nonumber\\
&\le& \sum_{n=0}^N\frac{\| D^n\varphi(Z_s)\|}{n!}\|\Delta_\star^{(n-1)}\|\,\Bigg\{\sum_{I\subsetneq \{1,\ldots,n\}}\hskip -3mm \vert  I \vert\, \bigg(\!\mop{max}\left(\big \| \pi_r ( \mathbf Z_t)\big\|,\big\|\pi_r(\mathbf Z_t^{'})\big\|\right) + C\left(T,\|\mathbf X\|_{\omega,\gamma}\right)\omega(s,t)^{\gamma}\bigg)^{\vert I\vert - 1}\nonumber\\
&&\hskip 5mm\times \big( \| \mathbf Z_t^{'} - \mathbf Z_t\| +\| \mathbf Z_s^{'} - \mathbf Z_s\| \big)\mop{max}\Big(\| R_{\mathbf X}(\mathbf Z)\|_{\omega, \gamma}, \| R_{\mathbf X}(\mathbf Z')\|_{\omega, \gamma}\Big)^{n- \vert I\vert}\hskip -2mm \omega(s,t)^{(N-r)\gamma(n-\vert I\vert - 1)}\nonumber\\
& & +\hskip -3mm\sum_{I\subsetneq \{1,\ldots,n\}} \hskip -3mm(n -\vert I\vert)\bigg(\!\mop{max}\left(\big \| \pi_r ( \mathbf Z_t)\big\|,\big\|\pi_r(\mathbf Z_t^{'})\big\|\right) + C\left(T,\|\mathbf X\|_{\omega,\gamma}\right)\omega(s,t)^{\gamma}\bigg)^{\vert I\vert}\nonumber\\
&&\hskip 5mm\times \|R_{\mathbf X}(\mathbf Z^{'}) - R_{\mathbf X}(\mathbf Z)\|_{\omega, \gamma} \ \mop{max}\Big(\| R_{\mathbf X}(\mathbf Z)\|_{\omega,\gamma}, \| R_{\mathbf X}(\mathbf Z')\|_{\omega,\gamma}\Big)^{n- \vert I\vert-1}\omega(s,t)^{(N-r)\gamma(n -\vert I\vert - 1)}\Bigg\}\nonumber\\
&&+\sum_{n=0}^{N} \frac{\|D^n\varphi(Z'_s)-D^n\varphi(Z_s)\|}{n!}\sum_{I\subsetneq\{1,\ldots,n\}}|I| \bigg(\mop{max}\left(\big \| \pi_r ( \mathbf Z_t)\big\|,\big\|\pi_r(\mathbf Z_t^{'})\big\|\right)+ C\left(T,\|\mathbf X\|_{\omega,\gamma}\right)\omega(s,t)^{\gamma}\bigg)^{|I|}\nonumber\\
&&\hskip 25 mm \times \mop{max}\Big(\left\|R_{\mathbf X}(\mathbf Z)\right\|_{\omega,\gamma},\,\left\|R_{\mathbf X}(\mathbf Z')\right\|_{\omega,\gamma}\Big)^{n-|I|}
\omega(s,t)^{(N-r)(n-|I|-1)\gamma}.
\end{eqnarray}
\ignore{
\noindent From the definition of the sup norm, we have for any controlled rough path $\mathbf Y$:

\begin{eqnarray}\label{est-three}
\left\| \mathbf Y\right\|_{\smop{sup}}&=&\sup_{t\in [0, T]}\left\| \mathbf Y_t\right\|\nonumber\\
& =& \sup_{t\in [0, T]}\sup_{h\neq 0} \frac{\left\| \langle h, \mathbf Y_t\rangle\right\|}{\left\| h \right\|}\nonumber \\
& = & \sup_{t\in [0, T]} \sup_{h\neq 0}  \frac{\left\| \langle h, R(\mathbf Y)_{0t}\rangle\right\|}{\left\| h \right\|}\nonumber\\
& \le & \sup_{s,t\in[0,T]}\sup_{h\neq 0} \frac{\left\| \langle h, R(\mathbf Y)_{st}\rangle\right\|}{\left\| h \right\|}\nonumber\\
& \le &\sup_{s,t\in[0,T]}\sup_{j}\left\| \pi_{j} R(\mathbf Y)_{st}\right\|\omega(s,t)^{-(N-j)\gamma}\omega(s,t)^{(N-j)\gamma}\nonumber\\
& \le &\left\| R(\mathbf Y)\right\|_{\omega, \gamma}\sup_{s,t\in[0,T]}\sup_{j}\omega(s,t)^{(N-j)\gamma}\nonumber\\
& \le &\left\| R(\mathbf Y)\right\|_{\omega, \gamma}\max\big(1, \omega(0,T)\big)^{N\gamma}\nonumber\\
& \le &K(T)\left\| \mathbf Y \right\|,
\end{eqnarray}\
where $K(T):=\max\big( 1, \omega(0,T)\big)^{\gamma}$.
}
\ignore{
We obviously have \textcolor{red}{!!!}
\begin{eqnarray*}
\|\mathbf Z\| = \|\mathbf Z_0\| + \|R_{\mathbf X}(\mathbf Z)\|_{\omega, \gamma}
 \le \|R_{\mathbf X}(\mathbf Z)\|_{\omega, \gamma},
\end{eqnarray*}
and similarly for $\mathbf Z'$. From $\mathbf Z_0=\mathbf Z'_0=0$ and $\varphi(\mathbf Z')_0=\varphi(\mathbf Z)_0$ \textcolor{red}{!!!}
}
We have
\[\| \mathbf Z^{'} - \mathbf Z\|=\|\mathbf Z'_0-\mathbf Z_0\|+\left\|R_{\mathbf X}(\mathbf Z')-R_{\mathbf X}(\mathbf Z)\right\|_{\omega,\gamma}\]
and
\[\|\varphi(\mathbf Z^{'}) - \varphi(\mathbf Z)\|=\|\varphi(\mathbf Z'_0)-\varphi(\mathbf Z_0)\|+\left\|R_{\mathbf X}\big(\varphi(\mathbf Z')\big)-R_{\mathbf X}\big(\varphi(\mathbf Z)\big)\right\|_{\omega,\gamma}.\]
Let us recall that, for any controlled rough path $\mathbf Y$, Estimate \eqref{sup-crp} reads
\begin{equation}\label{est-three}
\|\mathbf Y\|_{\smop{sup}}\le (N+1)K(T)^N\left(1+\|\mathbf X\|_{\omega,\gamma}\right)\|\mathbf Y\|
\end{equation}
with $K(T)=\mop{max}\big(1,\omega(0,T)\big)^{\gamma}$. From estimates \eqref{est-two} and \eqref{est-three}, and supposing we have chosen a graded norm on $\mathcal H^e$, we finally get the estimate
\allowdisplaybreaks{
{\small
\begin{eqnarray*}
&& \|\varphi(\mathbf Z')-\varphi(\mathbf Z)\|\le \sum_{n=0}^N\mopl{sup}_{s\in[0,T]}\frac{\| D^n\varphi(Z_s)\|}{n!}\|\Delta_\star^{(n-1)}\| \nonumber\\
&&\Bigg\{\sum_{I\subsetneq \{1,\ldots, n\}}\vert I\vert\Big(\max(\|\mathbf Z\|_{\smop{sup}},\|\mathbf Z'\|_{\smop{sup}})+C(T)K(T)\Big)^{\vert I\vert -1}(2\|\mathbf Z'-\mathbf Z\|_{\smop{sup}})K(T)^{N(n-\vert I\vert -1)}\mop{max}(\|\mathbf Z\|,\,\|\mathbf Z'\|)^{n-\vert I\vert}\nonumber\hskip -20mm\\
&&+\sum_{I\subsetneq \{1,\ldots, n\}}(n-\vert I\vert)\Big(\max(\|\mathbf Z\|_{\smop{sup}},\|\mathbf Z'\|_{\smop{sup}})+C\left(T,\|\mathbf X\|_{\omega,\gamma}\right)K(T)\Big)^{\vert I\vert }\|\mathbf Z'-\mathbf Z\|K(T)^{N(n-\vert I\vert -1)}\mop{max}(\|\mathbf Z\|,\,\|\mathbf Z'\|)^{n-\vert I\vert -1}\Bigg\}\nonumber\\
&&+\sum_{n=0}^{N} \frac{\|D^n\varphi(Z'_s)-D^n\varphi(Z_s)\|}{n!}\sum_{I\subsetneq\{1,\ldots,n\}}|I| \bigg(\mop{max}\left(\|\mathbf Z\|_{\smop{sup}},\,\|\mathbf Z'\|_{\smop{sup}}\right)+ C\left(T,\|\mathbf X\|_{\omega,\gamma}\right)K(T)^{\gamma}\bigg)^{|I|}\nonumber\\
&&\hskip 80 mm \times K(T)^{N(n-|I|-1}\mop{max}\big(\left\|\mathbf Z\right\|,\,\left\|\mathbf Z'\right\|)^{n-|I|}\\
&\le &\sum_{n=0}^N\mopl{sup}_{s\in[0,T]}\frac{\| D^n\varphi(Z_s)\|}{n!}\|\Delta_\star^{(n-1)}\|\nonumber\Bigg\{\sum_{I\subsetneq \{1,\ldots, n\}}\vert I\vert\left((N+1)(1+\|\mathbf X\|_{\omega,\gamma})K(T)^N(\|\mathbf Z\|+\|\mathbf Z'\|)+C\left(T,\|\mathbf X\|_{\omega,\gamma}\right)K(T)\right)^{\vert I\vert -1}\\
&&\hskip 65mm\times (2\|\mathbf Z'-\mathbf Z\|)(N-1)K(T)^NK(T)^{N(n-\vert I\vert -1)}(\|\mathbf Z\|+\|\mathbf Z'\|)^{n-\vert I\vert}\nonumber\\
&&+\sum_{I\subsetneq \{1,\ldots, n\}}(n-\vert I\vert)\left((N+1)(1+\|\mathbf X\|_{\omega,\gamma})K(T)^N(\|\mathbf Z\|+\|\mathbf Z'\|)+C\left(T,\|\mathbf X\|_{\omega,\gamma}\right)K(T)\right)^{\vert I\vert }\\
&&\hskip 80mm\times  \|\mathbf Z'-\mathbf Z\|K(T)^{N(n-\vert I\vert -1)}(\|\mathbf Z\|+\|\mathbf Z'\|)^{n-\vert I\vert -1}\Bigg\}\nonumber\\
&&+\sum_{n=0}^{N} \frac{\|D^n\varphi(Z'_s)-D^n\varphi(Z_s)\|}{n!}\sum_{I\subsetneq\{1,\ldots,n\}}|I| \bigg((N+1)K(T)^N\left(1+\|\mathbf X\|_{\omega,\gamma}\right)\left(\big  \|\mathbf Z\|+\|\mathbf Z'\|\right)+ C\left(T,\|\mathbf X\|_{\omega,\gamma}\right)K(T)\bigg)^{|I|}\\
&&\hskip 80mm \times K(T)^{N(n-|I|-1}\big(\left\|\mathbf Z\right\|+\left\|\mathbf Z'\right\|)^{n-|I|}\nonumber\\
&\le & CK(T)^{N^2}\|\varphi\|^N_{K(Z,Z')}\sum_{i,j,k=0}^{2N-1}\|\mathbf X\|_{\omega,\gamma}^i\|\mathbf Z\|^j\|\mathbf Z'\|^k\|\mathbf Z'-\mathbf Z\|  \hspace{2cm} \hbox{\ (from \eqref{degree-two})},
\end{eqnarray*}
}}
where $K(Z,Z'):=Z([0,T])\cup Z'([0,T])$, and where $\|\varphi\|^N_{K(Z,Z')}$ is given by \eqref{seminorms} with $K=K(Z,Z')$. This finally yields
\begin{equation*}
\|\varphi(\mathbf Z')-\varphi(\mathbf Z)\|\le C_{\varphi}(T)\|\mathbf Z'-\mathbf Z\|
\end{equation*}
with
\begin{equation}\label{Kphi}
C_\varphi(T):=CK(T)^{N^2}\|\varphi\|_{K(Z,Z')}^N\sum_{i,j,k=0}^{2N-1}\|\mathbf X\|_{\omega,\gamma}^i\|\mathbf Z\|^j\|\mathbf Z'\|^k.
\end{equation}

One can notice that the bounds $M(T)$, $K(T)$, $C(T)$, $C_q(T)$ and $C_{\varphi}(T)$ are nondecreasing with respect to the upper end $T$ of the chosen interval. This finishes up the proof of Proposition \ref{main-est-varphi}.
\end{proof}

\subsection{The sewing lemma and integration of controlled rough paths}\label{ssect:integration}
Let $\mathcal H$ be a connected graded commutative $N$-truncated Hopf algebra on the field $\mathbb R$ of real numbers, endowed with
\begin{itemize}
\item a nondegenerate homogeneous pairing $\langle -,-\rangle$ on $\mathcal H$.
\item a family of linear maps $L_i:\mathcal H\to \mathcal H$ indexed by $\{1,\ldots,d\}$, homogeneous of degree one, satisfying the cocycle condition
\[\Delta L_i(x)=(\mop{Id}\otimes L_i)\Delta(x)+L_i(x)\otimes \mathbf 1.\]
One moreover supposes that
\begin{itemize}
\item  The collection $\{L_i(\mathbf 1),\, i=1,\ldots,d\}$ is a basis of $\mathcal H_1$ (in coherence with Notation \ref{notation-sec-two}),
\item For any $i,j\in A
\{1,\ldots, d\}$ and $x,y\in\mathcal H$ we have
\begin{equation}\label{isom-Li}
\langle L_i(x),\, L_j(y)\rangle = \delta_i^j \langle x,y\rangle.
\end{equation}
In particular, the subspaces $T_i:=L_i(\mathcal H)$ are mutually orthogonal with respect to the pairing.
\end{itemize}
\end{itemize}
The dual $\mathcal H^\star$ will be identified to $\mathcal H$ via the linear isomorphism $x\mapsto \langle x,-\rangle$. We therefore have the direct orthogonal sum decomposition
\begin{equation}\label{orth}
\mathcal H=\mathcal H_0\oplus T\oplus F,\text{ where }\, T:=\bigoplus_{i\in A}^\perp T_i, \quad F:=T^\perp\cap\mop{Ker}\epsilon.
\end{equation}
In this context, the notation $\perp$ above the direct sum indicates that the summands are mutually orthogonal.\\
\begin{remark}\rm
The choice of notations is reminiscent to trees and forests respectively. In the case of the Butcher-Connes-Kreimer Hopf algebra $H=H_{\smop{BCK}}^{\{1,\ldots, d\}}$, with $\mathcal H=H_{\le N}$, the cocycle $L_i$ is the grafting operator $B_+^i$ on a common $i$-decorated root, the subspace $T_i$ is the linear span of trees with at most $N$ vertices whose root is decorated by $i$, and the subspace $F$ is the linear span of non-trivial forests, i.e. products of at least two (non-empty) rooted trees, with at most $N$ vertices. The same picture applies to the Munthe-Kaas--Wright Hopf algebra $H=H_{\smop{MKW}}^{\{1,\ldots, d\}}$ of $[d]$-decorated planar forests. The shuffle Hopf algebra $H=H_{\sqcup\hskip -1.5pt\sqcup}^{\{1,\ldots, d\}}$ endowed with the cocycles $L_i$ which add the letter $i$ on the right of a word also belongs to this framework: the third component $F$ is reduced to $\{0\}$ in this last example.
\end{remark}
Recall that $\mathbf X$ is an $N$-truncated rough path $\gamma$-adapted to $\omega$, and that $\mathbf Z$ is an $\mathbf X$- controlled rough path above $Z:[0,T]\to\mathbb R^e$.
\begin{lemma}\label{lem:deltaZ}
For any $i\in\{1,\ldots, d\}$, let $\Xi^i: \Delta_T\to \mathbb R^e$ be defined by
\[\Xi^i_{st}:= \langle \mathbf  X_{st},\, L_i(\mathbf  Z_s)\rangle = \langle ^t \!L_i(\mathbf  X_{st}),\, \mathbf  Z_s \rangle.\]
Then, for any $0\le s\le u\le t\le T$ we have
\begin{equation*}
\delta\Xi^i_{sut}:=\Xi^i_{su}+\Xi^i_{ut}-\Xi^i_{st} = \left\langle \mathbf X_{ut},\, L_i \big(R_{\mathbf X}(\mathbf Z)_{su}\big)\right\rangle.
\end{equation*}
\end{lemma}
\begin{proof} Using the cocycle condition for the $L_i$'s, we can compute
\begin{eqnarray*}
 \Xi^i_{su}+\Xi^i_{ut}-\Xi^i_{st}  &=&  \langle \mathbf  X_{su}, L_{i}(\mathbf  Z_{s}) \rangle + \langle \mathbf  X_{ut}, L_{i}(\mathbf  Z_{u})\rangle-\langle \mathbf  X_{st}, L_{i}(\mathbf  Z_{s}) \rangle  \\
&=&   \langle \mathbf  X_{su}, L_{i}(\mathbf  Z_{s})\rangle + \langle \mathbf  X_{ut}, L_{i}(\mathbf  Z_{u})\rangle -\langle \mathbf  X_{su} \star \mathbf  X_{ut}, L_{i}(\mathbf  Z_{s})\rangle \\
&=&   \langle \mathbf  X_{su}, L_{i}(\mathbf  Z_{s})\rangle + \langle \mathbf  X_{ut}, L_{i}(\mathbf  Z_{u})\rangle-\langle \mathbf  X_{su} \otimes \mathbf  X_{ut}, \Delta L_{i}(\mathbf  Z_{s})\rangle \\
&=&  \langle \mathbf  X_{su}, L_{i}(\mathbf  Z_{s})\rangle + \langle \mathbf  X_{ut}, L_{i}(\mathbf  Z_{u})\rangle-\langle \mathbf  X_{su} \otimes \mathbf  X_{ut}, (I \otimes L_{i})\Delta \mathbf  Z_{s} + L_{i}(\mathbf  Z_{s}) \otimes \mathbf 1\rangle    \\
&=& \langle \mathbf  X_{ut}, L_{i}(\mathbf  Z_{u})\rangle  -\langle \mathbf  X_{su} \star \,  ^t \! L_{i}(\mathbf  X_{ut}), \mathbf  Z_{s}\rangle   \\
&=&  \langle ^t\!L_{i}(\mathbf  X_{ut}), T_{su} \mathbf Z_s -\mathbf  Z_{u}\rangle\\
&=& \left\langle \mathbf  X_{ut},\, L_i\big(R_{\mathbf X}(\mathbf Z)_{su}\big)\right\rangle.   
\end{eqnarray*} 
\end{proof}

\noindent The map $(s,t)\mapsto \Xi_{st}^i$ is \textsl{almost} an increment. More precisely,

\begin{theorem}\label{delta-xi}
\begin{equation}\label{gen-sewing}
\| \delta\Xi_{sut}^i\| \le \|L_i\|\,\|\mathbf X\|_{\omega,\gamma}\,\|\mathbf Z\|\sum_{j=1}^N\omega(u,t)^{\alpha_j}\omega(s,u)^{\beta_j}
\end{equation}
with $\alpha_j+\beta_j=1+\varepsilon$ with $\varepsilon>0$.
\end{theorem}
\begin{proof}
This estimate is proven by a direct computation:
\begin{eqnarray*}
\|\delta\Xi^i_{sut}\|=\| \Xi_{su}^i+\Xi_{ut}^i-\Xi_{st}^i \| &=&\| \langle \mathbf  X_{ut},\,L_i\big(R_{\mathbf X}(\mathbf Z)_{su}\big)\rangle \|\\
&\le& \sum_{j=1}^N\| \langle \pi_j^*(\mathbf  X_{ut}),\,\pi_jL_i\big(R_{\mathbf X}(\mathbf Z)_{su}\big)\rangle \|\\
&\le& \sum_{j=1}^N\| \langle \pi_j^*(\mathbf  X_{ut}),\,L_i\pi_{j-1}\big(R_{\mathbf X}(\mathbf Z)_{su}\big)\rangle \|\\
&\le &\|L_i\|\,\|\mathbf X\|_{\omega,\gamma}\|\,\|R_{\mathbf X}(\mathbf Z)\|_{\omega,\gamma}\sum_{j=1}^N\omega(u,t)^{j\gamma}\omega(s,u)^{(N+1-j)\gamma}\\
&\le &\|L_i\|\,\|\mathbf X\|_{\omega,\gamma}\|\,\|\mathbf Z\|\sum_{j=1}^N\omega(u,t)^{j\gamma}\omega(s,u)^{(N+1-j)\gamma}.
\end{eqnarray*}
This gives \eqref{gen-sewing}, with $\varepsilon=(N+1)\gamma-1$.
\end{proof}
\begin{remark}\rm
The operators $L_i$ increase the degree by one, and their tranpose $^t \!L_i$ lower the degree by one. In particular, $L_i(\mathbf Z_s)$ does not depend on the $N$-th homogeneous component $(\mathbf Z_s)_N$ (the image of it by $L_i$ is of degree $N+1$, hence $=0$). We can therefore safely suppose $(\mathbf Z_s)_N=0$, which explains why controlled rough paths are $(N-1)$-truncated in the literature (see e.g. \cite[p.6270, Definition 4.10]{FZ18} and~\cite[p.79, Definition 3.3.3]{Kel}).
\end{remark}

Let us recall the generalized sewing lemma \cite[Theorem 2.5]{FZ18}, where the Refinement Riemann-Stiejes (RRS) convergence of a Riemann sum to a limit $K$ is defined as follows: for any $\varepsilon > 0$, there exists a partition $P_{\varepsilon}$ such that for any refinement $P\supset P_{\varepsilon}$, one has $$\left\|\left(\sum_{[u,v]\in P}\Xi_{u,v}\right)-K\right\|< \varepsilon,$$
where $\Xi$ is a map from $\Delta_T$ into a normed vector space $(V,\|-\|)$. The element $K\in V$ is the RRS-limit of the Riemann sums of $\Xi$.

\begin{theorem}\label{thm:gen-sewing}
(Generalized sewing lemma) Let $\Xi$ be a map from $\Delta_T$ into $\mathbb R^e$ with
\begin{align*}
\| \Xi_{st} - \Xi_{su} - \Xi_{ut} \|\leq \sum_{j=1}^{n}\omega_{1,j}^{\alpha_1, j}(s,u)\omega_{2,j}^{\alpha_{2},j}(u,t),
\end{align*}
where $\omega_{1,j}, \omega_{2,j}$ are controls and $\alpha_{1,j} + \alpha_{2,j}> 1$ for all $j=1,\ldots, n$. Then the following (unique) limit exists:
\begin{align*}
\mathcal{I}\Xi_{st}:= \mopl{RRS-lim}_{P\smop{partition of }[s,t]}\sum_{(u,v)\in P} \Xi_{uv}.
\end{align*}
This limit is additive in the sense that $\mathcal{I}\Xi_{st}=\mathcal{I}\Xi_{su}+\mathcal{I}\Xi_{ut}$ for any $0\le s\le u\le t\le T$,
and one has the following estimate:
\[\| \mathcal{I}\Xi_{st} - \Xi_{st}\|\leq C\sum_{j=1}^{n}\omega_{1,j}^{\alpha_{1},j}(s,t-)\omega_{2,j}^{\alpha_{2},j}(s+,t),\]
with $C$ depending only on $\min\{\alpha_{1,j} + \alpha_{2,j}, j=1, \ldots, n\}$.
\end{theorem}

From the generalized sewing lemma (Theorem \ref{thm:gen-sewing}), returning to the notations of Lemma \ref{lem:deltaZ}, for any $i\in\{1,\ldots, d\}$, there exists a remainder $r^i_{st}$ with
\[
\vert r_{st}^i\vert\le nC\omega(s,t)^{1+\varepsilon}
\]
such that $\mathcal I\Xi_{st}^i:= \Xi_{st}^i+r_{st}^i$ is an increment, \textsl{the rough integral from $s$ to $t$ of $Z$ against $X^i$}, denoted by\footnote{The notation is somewhat misleading, as this integral does depend on the choice of a rough path $\mathbf X$ above $X$.}
\[
\mathcal I Z_{st}^i=:\int_s^t Z_r\,dX_r^i=:\int_0^tZ_r\, dX_r^i-\int_0^sZ_r\, dX_r^i,
\]
The estimate \eqref{gen-sewing} permits us to apply the generalized sewing lemma (Theorem \ref{thm:gen-sewing}), with all controls $\omega_{1,j},\omega_{2,j}$ equal to the same control $\omega$.
The next step is to build up an $\mathbf  X$-controlled rough path
$$t\mapsto \int_0^t \mathbf  Z_r\, dX_r^i\,\text{ above }\, t\mapsto \int_0^t Z_r\, dX_r^i.$$ The definition is given by
\begin{equation*}
\int_0^t\mathbf Z_r\, dX_r^i:=L_i(\mathbf Z_t)+\left(\int_0^t Z_r\,dX_r^i\right)\mathbf 1
\end{equation*}
for any $t\in[0,T]$. We therefore have in view of \eqref{isom-Li} and \eqref{orth}:
\begin{eqnarray*}
\left\langle \mathbf 1,\, \int_0^t \mathbf  Z_r\, dX_r^i\right\rangle &=& \int_0^t Z_r\, dX_r^i \hbox{ for any } i\in \{1,\ldots,d\}, \label{above-Z}\\
\left\langle L_j(x),\, \int_0^t \mathbf  Z_r\, dX_r^i\right\rangle &=& \delta_i^j\langle x,\, \mathbf  Z_t \rangle \hbox{ for any } i,j\in \{1,\ldots,d\},\\
\left\langle y,\, \int_0^t \mathbf  Z_r\, dX_r^i\right\rangle &=& 0 \hbox { for any } y\in F \hbox { and } i\in \{1,\ldots,d\}.
\end{eqnarray*}

\noindent We can now state the main result of this section.

\begin{proposition}\label{crpint}
The $\mathcal H^e$-valued path $\int_0^\bullet \mathbf Z_r\, dX_r^i$ is an $\mathbf X$-controlled rough path above $\int_0^\bullet  Z_r\, dX_r^i$.
\end{proposition}

\begin{proof}
We check \eqref{crp} by direct computation:
\begin{eqnarray*}
\left\langle \mathbf X_{st}\star h,\,\int_0^s \mathbf Z_r\, dX_r^i\right\rangle &=& \left\langle\mathbf X_{st}\star h,\, L_i(\mathbf Z_s)+\left(\int_0^s Z_r\, dX_r^i\right)\mathbf 1\right\rangle\\
&=&\left\langle\mathbf X_{st}\otimes h,\, \Delta\big(L_i(\mathbf Z)\big)+\left(\int_0^s Z_r\, dX_r^i\right)\mathbf 1\otimes\mathbf 1\right\rangle \\
&=&\left(\int_0^s Z_r\, dX_r^i\right)\langle h,\,\mathbf 1\rangle+\langle \mathbf X_{st}\otimes h,\, (\mop{Id}\otimes L_i)\Delta(\mathbf Z_s)+L_i(\mathbf Z_s)\otimes\mathbf 1\rangle\\
&=&\left(\int_0^s Z_r\, dX_r^i\right)\langle h,\,\mathbf 1\rangle+\langle \mathbf X_{st}\otimes h,\, (\mop{Id}\otimes L_i)\Delta(\mathbf Z_s)\rangle+\left(\int_s^t Z_r\, dX^i- r^i_{st}\right)\langle h,\mathbf 1\rangle \\
&=&\left(\int_0^t Z_r\, dX_r^i- r^i_{st}\right)\langle h,\mathbf 1\rangle + \langle \mathbf X_{st}\otimes h,\, (\mop{Id}\otimes L_i)\Delta(\mathbf Z_s)\rangle.\\
\end{eqnarray*}
In view of \eqref{orth}, we therefore get
\[\left\langle \mathbf X_{st},\,\int_0^s \mathbf Z_r\, dX_r^i\right\rangle=\int_0^t Z_r\, dX_r^i-r_{st}^i=\left\langle \mathbf 1,\,\int_0^t \mathbf Z_r\, dX_r^i\right\rangle-r_{st}^i\]
and
\[\left\langle \mathbf X_{st}\star h,\,\int_0^s \mathbf Z_r\, dX_r^i\right\rangle=0=\left\langle h,\,\int_0^t \mathbf Z_r\, dX_r^i\right\rangle\]
for any $h\in F$ or $h\in L_j(\mathcal H)$ with $j\neq i$. For $h=L_i(y)$ we get
\begin{eqnarray}
\left\langle \mathbf X_{st}\star h,\,\int_0^s \mathbf Z_r\, dX_r^i\right\rangle&=&\langle \mathbf X_{st}\otimes L_i(y),\, (\mop{Id}\otimes L_i)\Delta(\mathbf Z_s)\rangle\nonumber\\
&=&\langle \mathbf X_{st}\otimes y,\, \Delta(\mathbf Z_s)\rangle\nonumber\\
&=&\langle\mathbf X_{st}\star y,\, \mathbf Z_s\rangle\nonumber\\
&=&\langle y,\, \mathbf Z_t-R_{\mathbf X}(\mathbf Z)_{st}\rangle\nonumber\\
&=&\left\langle h,\int_0^t \mathbf Z_r\, dX_r^i\right\rangle -\langle y,\,R_{\mathbf X}(\mathbf Z)_{st}\rangle.\label{remainder-int}
\end{eqnarray}
We have
\[\|\langle y,\,R_{\mathbf X}(\mathbf Z)_{st}\rangle\|_{\mathbb R^e}\le C\omega(s,t)^{(N-|y|)\gamma}\le  C\omega(s,t)^{(N-|h|)\gamma},\]
which ends up the proof of Proposition \ref{crpint}.
\end{proof}
\begin{remark}
Identity \eqref{remainder-int} can be re-phrased as
\begin{equation}\label{remainder-int-bis}
R_{\mathbf X}\left(\int_0^\bullet \mathbf Z_r\,dX_r^i\right)_{st}=r_{st}^i\mathbf 1+L_i\big(R_{\mathbf X}(\mathbf Z)_{st}\big).
\end{equation}
\end{remark}
\begin{theorem}\label{ixi-xi}
Keeping the previous notations, the following estimate holds:
\begin{equation*}
\| \mathcal I\Xi_{st}^i-\Xi_{st}^i\| \le C\|L_i\|\,\|\mathbf X\|_{\omega,\gamma}\,\|\mathbf Z\|\sum_{j=1}^N\omega(s,t_-)^{j\gamma}\omega(s_+,t)^{(N+1-j)\gamma},
\end{equation*}
where $C$ is a combinatorial constant.
\end{theorem}
\begin{proof}
It is a direct consequence of Theorem \ref{delta-xi} applied to $\Xi^i$, and Theorem \ref{thm:gen-sewing} applied to $\Xi^i/(\|L_i\|\,\,\|\mathbf X\|_{\omega,\gamma}\,\|\mathbf Z\|)$, in the particular case when all controls $\omega_{1,j}$ and $\omega_{2,j}$ coincide with $\omega$.
\end{proof}
\begin{corollary}\label{ixi-xi-bis}
With the preceding notations,
\[\| \mathcal I\Xi_{st}^i-\Xi_{st}^i\| \le NC\|L_i\|\,\|\mathbf X\|_{\omega,\gamma}\,\|\mathbf Z\|\omega(s,t)^{(N+1)\gamma}.\]
\end{corollary}
\noindent We end up this section with a lemma for later use in Section \ref{ult-bis}.
\begin{lemma}\label{Li-omega-gamma}
Let $R:\Delta_T\to \mathcal H^e$. Then for any $i=1,\ldots, d$ we have
\[\|L_i\circ R\|_{\omega,\gamma}\le \omega(0,T)^\gamma\|L_i\|\,\|R\|_{\omega,\gamma}.\]
\end{lemma}
\begin{proof}
Given by the straightforward computation below, using $\pi_j\circ L_i=L_i\circ \pi_{j-1}$ for any $j=1,\ldots,N$, due to the fact that $L_i$ raises the degree by one:
\begin{eqnarray*}
\|L_i\circ R\|_{\omega,\gamma}&=&\mopl{sup}_{j=1,\ldots, N}\,\mopl{sup}_{s,t\in\Delta_T}\left\|\pi_j\circ L_i(R)_{st}\right\|\omega(s,t)^{-(N-j)\gamma}\\
&=&\mopl{sup}_{j=1,\ldots, N}\,\mopl{sup}_{s,t\in\Delta_T}\left\|L_i\circ \pi_{j-1}(R)_{st}\right\|\omega(s,t)^{-(N-j)\gamma}\\
&\le &\|L_i\|\mopl{sup}_{j=1,\ldots, N}\,\mopl{sup}_{s,t\in\Delta_T}\left\|\pi_{j-1}(R)_{st}\right\|\omega(s,t)^{-(N-j-1)\gamma}\omega(s,t)^{\gamma}\\
&\le& \omega(0,T)^\gamma\|L_i\|\,\|R\|_{\omega,\gamma}.
\end{eqnarray*}
\end{proof}
\section{Existence and uniqueness for the initial value problem}\label{sect:ult}
In this section, we prove existence and uniqueness of the solution of the initial value problem \eqref{ivplifted} for truncated rough path with jumps and any roughness,
based on the framework of truncated Hopf algebras. Let $N$ be a positive integer, let $\mathcal H$ be an $N$-truncated connected graded commutative Hopf algebra fulfilling the properties described in Paragraph \ref{ssect:integration}. Let $d$ be the dimension of $\mathcal H_1$, let $T>0$, and let $\omega:\Delta_{T}\to [0,+\infty[$ be a right-continuous control, in the sense that
\[\omega(s,t_+)=\omega(s,t)\]
for any $0\le s<t\le T$. Let $\gamma>0$ and let $\mathbf X$ be an $N$-truncated rough path  $\gamma$-adapted to the control $\omega$.\\

\noindent Given the preceding setup, we proceed to present two key results of the present paper.
\begin{theorem}\label{thm:local}(Local existence and uniqueness)
For any positive sufficiently small $T'$, the initial value problem \eqref{ivplifted} admits a unique solution in the controlled rough path space $\mathcal V_{\mathbf X}^{[0,T']}$.
\end{theorem}
\begin{proof}
\noindent We denote the integration operators by $\mathbf Z\mapsto \mathcal I^i_{\mathbf X}(\mathbf Z)$, with
\begin{eqnarray*}
\mathcal I^i_{\mathbf X}(\mathbf Z)_t := \int^t_0\,\mathbf Z_r dX_r^i= L_i(\mathbf Z_t)+\left(\int_0^t Z_r\,dX_r^i\right)\mathbf 1.
\end{eqnarray*}

From Proposition \ref{crpint}, the integration operator $\mathcal I^i_{\mathbf X}$ is a linear operator on the controlled rough paths space $\mathcal V_{\mathbf X}^{[0,T']}$. Recall from (\ref{crp-norm}) that the norm of controlled rough path with hypothesis $\mathbf Y_0=0$ is:
\begin{eqnarray*}
 \| \mathcal I^i_{\mathbf X}(\mathbf Y)\| =  \|R\big(\mathcal I^i_{\mathbf X}(\mathbf Y)\big)\|_{\omega, \gamma}
  =  \mopl{sup}_{(s,t)\in \Delta T}\mopl{sup}_{j=0,\ldots, N}\|\pi_j \circ R\big(\mathcal I^i_{\mathbf X} (\mathbf Y)\big)_{s,t}\|_{\omega, \gamma}.
\end{eqnarray*}
From the proof in proposition \ref{crpint}, we have
\begin{eqnarray*}
  \left\langle \mathbf X_{st}\star L_i(y), \int_0^s \mathbf Y_r\, dX_r^i \right\rangle &=& \left\langle h,\int^t \mathbf Y_r\, dX_r^i\right\rangle -R_{st}^y,
\end{eqnarray*}
for $L_i(y)=h$ and $R_{st}^y=\langle y, R(\mathbf Y)_{st}\rangle$. In view of (\ref{crp}), we also get
\begin{eqnarray*}
  \left\langle \mathbf X_{st}\star h, \int_0^s \mathbf Y_r\, dX_r^i\right\rangle &=& \left\langle h, \int_0^t \mathbf Y_r\, dX_r^i\right\rangle - \left\langle h, R\big( \mathcal I^i(\mathbf Y)\big)_{st}\right\rangle.
\end{eqnarray*}
Therefore
\begin{eqnarray*}
  \left\langle L_i(y), R\big(\mathcal I^i_{\mathbf X}(\mathbf Y)_{st}\big)\right\rangle &=& \left\langle y, R\big( \mathbf Y\big)_{st}\right\rangle.
\end{eqnarray*}
For $h\in F$ or $h=L_{i'}(y), i'\neq i$, we have $\|\langle h, R\big( \mathcal I^i_{\mathbf X}(\mathbf Y)_{st}\big)\rangle\| = 0$ and for $h=L_i(y)$, and so
\begin{eqnarray*}
  \left\langle L_i(y), R\big(\mathcal I^i_{\mathbf X}(\mathbf Y)_{st}\big)\right\rangle = \left\langle y, R\big( \mathbf Y\big)_{st}\right\rangle
 \le  \|R(\mathbf Y)\|_{\omega, \gamma}\omega(s,t)^{(N-j)\gamma},
\end{eqnarray*}
where we suppose $y\in\mathcal H_j$. Therefore
\[
  \|\pi_{j+1}\circ R\big(\mathcal I^i_{\mathbf X}(\mathbf Y)_{st}\big)\|
   \le  \|R(\mathbf Y)\|_{\omega, \gamma}\omega(s,t)^{(N-j)\gamma},
\]
which means
\begin{eqnarray*}
  \|\mathcal I^i_{\mathbf X}(\mathbf Y)\| = \|R\big(\mathcal I^i_{\mathbf X}(\mathbf Y)\big)\|_{\omega, \gamma} &\le& \|R(\mathbf Y)\|_{\omega, \gamma}\mopl{sup}_{s,t}\omega(s,t)^{\gamma}\\
  &\le& \|R(\mathbf Y)\|_{\omega, \gamma}\omega(0,T')^\gamma.
\end{eqnarray*}
We therefore get the estimate
\begin{eqnarray*}
  \|\mathcal I^i_{\mathbf X}(\mathbf Y)\| & \le & \omega(0, T')^\gamma \|\mathbf Y\|,
\end{eqnarray*}
with the hypothesis $\mathbf Y_0=0$. Choosing two controlled rough paths $\mathbf Z$ and $\mathbf Z'$ such that $\mathbf Z_0=\mathbf Z'_0$ we therefore get
\begin{equation}\label{intn}
 \|\mathcal I^i_{\mathbf X}(\mathbf Z'-\mathbf Z)\|  \le  \omega(0, T')^\gamma \|\mathbf Z'-\mathbf Z\|.
\end{equation}

Now let us choose two controlled rough paths $\mathbf Z,\mathbf Z'$ in a ball of radius $M$ (for the controlled rough path norm \eqref{crp-norm}) such that $\mathbf Z_0=\mathbf Z'_0$. We obviously have
$$(\mathcal M_{T'}\mathbf Z)_0=(\mathcal M_{T'}\mathbf Z')_0=\mathbf Z_0=\mathbf Z'_0,$$
where $\mathcal M_{T'}$ is the transformation of $\mathcal V_{\mathbf X}^{[0,T']}$ defined by \eqref{eq:firstm}.
Using the norm of integration (\ref{intn}) and Proposition \ref{main-est-varphi} for a fixed $t\in [0, T']$, we can get the estimate of $\|\mathcal M_{T'}\mathbf Z' - \mathcal M_{T'} \mathbf Z \|$ as follows
\begin{eqnarray}
  \| \mathcal M_{T'} \mathbf Z' - \mathcal M_{T'} \mathbf Z\| & = & \left\| \mathbf Z'_0 + \sum_{i=1}^d\int_0^t \varphi_i(\mathbf Z')_r \, d\mathbf X_r^i - \mathbf Z_0 - \sum_{i=1}^d\int_0^t \varphi_i(\mathbf Z)\, d\mathbf X_r^i\right\|\nonumber\\
  & = & \left\| \sum_{i=1}^d\int_0^t \varphi_i(\mathbf Z')_r \, d\mathbf X_r^i  - \sum_{i=1}^d\int_0^t \varphi_i(\mathbf Z)\, d\mathbf X_r^i \right\|\nonumber\\
  & = & \left\| \sum_{i=1}^d \int_0^t \Big(\varphi_i(\mathbf Z')_r - \varphi_i(\mathbf Z)_r\Big)\, d\mathbf X_r^i \right\|\nonumber\\
  & \le & \omega(0, T')^\gamma \left\| \sum_{i=1}^d \varphi_i(\mathbf Z') - \varphi_i(\mathbf Z)\right\|\nonumber\\
  & \le &\omega(0, T')^\gamma \sum_{i=1}^d C_{M, \varphi_i}(T')\| \mathbf Z' - \mathbf Z\|,\label{contraction}
\end{eqnarray}

where the dependence of the constants $C_{\varphi_i}(T')$ with respect to the radius $M$ is emphasized in the notation. Using the right-continuity of the control $\omega$ (more precisely, of the map $t\mapsto \omega(0,t)$ at $t=0$), one can see that $\mathcal M_{T'}$ is contracting provided the upper bound $T'$ of the interval is small enough. This finishes up the proof of the existence and uniqueness of a local solution of the lifted initial value problem \eqref{ivplifted}.
\end{proof}

\begin{theorem}(Global existence and uniqueness)\label{thm:ult}
The initial value problem \eqref{ivplifted} admits a unique solution in the controlled rough path space $\mathcal V_{\mathbf X}^{[0,T]}$ whenever the maps $\varphi_i$ and their derivatives up to order $N$ are bounded on $\mathbb R^e$.
\end{theorem}
\begin{proof}
We use the last item of Proposition \ref{prop-control}, and the local existence and uniqueness Theorem \ref{thm:local}. From the boudedness of the iterated derivatives of the $\varphi_i$'s, the bound $T'$ in the local existence and uniqueness theorem can be chosen independently from the initial condition\footnote{This boundedness condition is crucial even in the ordinary smooth case, as illustrated by the initial value problem $\dot y(t)=y(t)^2$ with initial condition $y(0)=c$: the unique solution $y(t)=c/(1 -ct)$ explodes at $t=1/c$.}. We choose $\delta>0$ such that
\[\delta^\gamma\sum_{i=1}^d C_{M,\varphi_i}(T') <1.\]

From right-continuity of the control, from Proposition \ref{prop-control} and from \eqref{contraction}, there exists a partition $\mathcal P=\{t_0,\ldots, t_k\}$ of the interval $[0,T]$ with
\[0=t_0<t_1<\cdots <t_k=T\]
such that $\omega(t_j,t_{j+1,-})<\delta$ for any $j=0,\ldots, k-1$. The initial value problem \eqref{ivplifted} therefore admits a unique solution $\mathbf Z$ on $[0,T']$ for any $T'<t_1$, hence on the first interval $[0,t_1[$
of the partition. Now the solution is obviously continuous, and \eqref{ivplifted} is still verified at $t_1$. Now we solve the initial value problem on $[t_1,t_2[$ with initial condition $\mathbf Z_{t_1}$, and extend it to $[t_1,t_2]$, again by continuity. We now get a unique solution on $[0,t_2]$. Proceeding the same way along the finite partition $\mathcal P$, we finally get the existence and uniqueness of the solution $\mathbf Z$ of \eqref{ivplifted} on the whole interval $[0,T]$.
\end{proof}
\begin{remark}
The degree zero part of the unique solution $\mathbf Z$ of \eqref{ivplifted} given by Theorem \ref{thm:ult} gives back a path $Z:[0,T]\to\mathbf R^e$, which deserves to be called a solution of the initial value problem \eqref{ivp}. The set of solutions of \eqref{ivp} is therefore parameterised by the set of rough paths $\mathbf X$ above the driving path $X$.
\end{remark}
\section{Universal Limit Theorem}\label{ult-bis}
We finally address the robustness of the solution of \eqref{ivp} with respect to the initial condition, the vector fields $(\varphi_i)_{i=1,\ldots, d}$ and the driving rough path $\mathbf X$.
\subsection{Formal distance and supremum norm}
Let $\mathbf Z$ and $\widetilde{\mathbf Z}$ be two controlled rough paths driven by two possibly different rough paths $\mathbf X$ and $\widetilde{\mathbf X}$ respectively. They belong to two different Banach spaces $\mathcal V_{\mathbf X}^{[0,T]}$ and $\mathcal V_{\widetilde{\mathbf X}}^{[0,T]}$ respectively, and their remainders are defined by
\[R_{\mathbf X}(\mathbf Z)_{st}:=\mathbf Z_t-T_{st}\mathbf Z_s,\hskip 12mm R_{\widetilde{\mathbf X}}(\widetilde{\mathbf Z})_{st}:=\widetilde {\mathbf Z_t}-\widetilde{T_{st}}\mathbf Z_s,\]
where $T_{st}$ is the componentwise application in $\mathcal H^e$ of the transpose $^t(\mathbf X_{st}\star -)$, and similarly for $\widetilde {T_{st}}$. The formal distance between both controlled paths \cite[Paragraph 4.3, p.6271]{FZ18} is defined by
\[d(\mathbf Z,\widetilde{\mathbf Z}):=\|\mathbf Z_0-\widetilde{\mathbf Z}_0\|_{\mathcal H^e}+\|R_{\mathbf X}(\mathbf Z)-R_{\widetilde{\mathbf X}}(\widetilde{\mathbf Z})\|_{\omega,\gamma}.\]
It coincides with the distance inside  $\mathcal V_{\mathbf X}^{[0,T]}$ whenever $\widetilde{\mathbf X}=\mathbf X$. The following proposition addresses the effect of a smooth map on this distance.
\begin{proposition}\label{d-varphi}
Let $\mathbf Y$ and $\widetilde{\mathbf Y}$ be two controlled rough paths on $[0,T]$ driven by two possibly different rough paths $\mathbf X$ and $\widetilde{\mathbf X}$ respectively, and let $\varphi\in C^N(\mathbb R^e,\,\mathbb R^e)$. Then we have
\[d\left(\varphi(\mathbf Y),\, \varphi(\widetilde{\mathbf Y})\right)\le C_{\varphi}(T)\,d\left(\mathbf Y,\,\widetilde{\mathbf Y}\right),\]
where the constant $C_{\varphi}(T)$ can be chosen non-decreasing with respect to $T$. It depends on $\|\mathbf X\|_{\omega,\gamma}$, $\|\widetilde{\mathbf X}\|_{\omega,\gamma}$, $\|\mathbf Y\|$, $\|\widetilde{\mathbf Y}\|$ and $\|\varphi\|_{K(Y,Y')}^N$.
\end{proposition}
\begin{proof}
This proposition generalises Proposition \ref {main-est-varphi} to the case when both controlled paths are driven by two different rough paths. Its proof is very similar to the proof of Proposition \ref{main-est-varphi} modulo the change of notations $\mathbf Z\mapsto \mathbf Y$ and $\mathbf Z'\mapsto \widetilde{\mathbf Y}$, except that the formal distance $d\left(\varphi(\mathbf Y),\, \varphi(\widetilde{\mathbf Y})\right)$ replaces the norm $\left\|\varphi(\widetilde{\mathbf Y})-\varphi(\mathbf Y)\right\|$. Details are left to the reader.
\end{proof}
\subsection{The universal limit theorem}
Now we suppose
\[\mathbf Z=\mathbf Y_0+\sum_{i=1}^d \mathcal I^i_{\mathbf X}\big(\varphi_i(\mathbf Y)\big),\hskip 12mm \widetilde{\mathbf Z}=\widetilde{\mathbf Y}_0+\sum_{i=1}^d \mathcal I^i_{\widetilde{\mathbf X}}\big({\widetilde\varphi}_i(\widetilde{\mathbf Y})\big),\]
which yields
\begin{equation*}
d(\mathbf Z,\widetilde{\mathbf Z})\le\|\mathbf Y_0-\widetilde{\mathbf Y}_0\|+\sum_{i=1}^d \left\|R_{\mathbf X}\Big(\mathcal I^i_{\mathbf X}\big(\varphi_i(\mathbf Y)\big)\Big)-R_{\widetilde{\mathbf X}}\Big(\mathcal I^i_{\widetilde{\mathbf X}}\big(\widetilde{\varphi}_i(\widetilde{\mathbf Y})\big)\Big)\right \|_{\omega,\gamma}.
\end{equation*}
In view of \eqref{remainder-int-bis}, we get, with
\[
r^i_{st}:=\int_s^t(\varphi_i\circ Y)_u\, dX_u^i-\left\langle \mathbf  X_{st},\, L_i\big(\varphi_i(\mathbf  Y_s)\big)\right\rangle, \hskip 12mm \widetilde r\,^i_{st}:=\int_s^t(\widetilde \varphi_i\circ \widetilde Y)_u\, d\widetilde X_u^i-\left\langle \widetilde {\mathbf  X}_{st},\, L_i\big(\widetilde \varphi_i(\widetilde {\mathbf  Y}_s)\big)\right\rangle,
\]
\allowdisplaybreaks{
\begin{eqnarray}
d(\mathbf Z,\widetilde{\mathbf Z})&\le&\|\mathbf Y_0-\widetilde{\mathbf Y}_0\|+\sum_{i=1}^d \left\|L_i\Big(R_{\mathbf X}\big(\varphi_i(\mathbf Y)\big)-R_{\widetilde{\mathbf X}}\big(\widetilde{\varphi}_i(\widetilde{\mathbf Y})\big)\Big)+(r^i-\widetilde r\,^{i})\mathbf 1 \right \|_{\omega,\gamma}\nonumber\\
&\le &\|\mathbf Y_0-\widetilde{\mathbf Y}_0\|+\sum_{i=1}^d \left\|L_i\Big(R_{\mathbf X}\big(\varphi_i(\mathbf Y)\big)-R_{\mathbf X}\big(\widetilde{\varphi}_i(\mathbf Y)\big)\Big)\right \|_{\omega,\gamma}+\sum_{i=1}^d \left\|L_i\Big(R_{\mathbf X}\big(\widetilde{\varphi}_i(\mathbf Y)\big)-R_{\widetilde{\mathbf X}}\big(\widetilde{\varphi}_i(\widetilde{\mathbf Y})\big)\Big)\right \|_{\omega,\gamma}\nonumber\\
&&\hskip 20mm + \sum_{i=1}^d \left\|(r^i-\widetilde r\,^{i})\mathbf 1 \right \|_{\omega,\gamma}\nonumber\\
&\le& \|\mathbf Y_0-\widetilde{\mathbf Y}_0\|+\omega(0,T)^\gamma\sum_{i=1}^d\|L_i\|\,\left\|R_{\mathbf X}\big((\varphi_i-\widetilde\varphi\, _i)(\mathbf Y)\big)\right\|_{\omega,\gamma}\nonumber\\
&&\hskip 20mm+\omega(0,T)^\gamma\sum_{i=1}^d\|L_i\|\,\left\|R_{\mathbf X}\big(\widetilde{\varphi}_i(\mathbf Y)\big)-R_{\widetilde{\mathbf X}}\big(\widetilde{\varphi}_i(\widetilde{\mathbf Y})\big)\right\|_{\omega,\gamma}\nonumber\\
&&\hskip 20mm + \sum_{i=1}^d \left\|(r^i-\widetilde r\,^{i})\mathbf 1 \right \|_{\omega,\gamma}\nonumber \hspace{2cm}\hbox{ (from Lemma \ref{Li-omega-gamma})}\\
&\le& \|\mathbf Y_0-\widetilde{\mathbf Y}_0\|+\omega(0,T)^\gamma\sum_{i=1}^d \|L_i\|\bigg(\left\|(\varphi_i-\widetilde\varphi_i)(\mathbf Y)\right\|+d\left(\widetilde\varphi_i(\mathbf Y),\,\widetilde\varphi_i(\widetilde{\mathbf Y})\right)\bigg)\nonumber\\
&&\hskip 20mm +\omega(0,T)^\gamma\sum_{i=1}^d \mopl{sup}_{s,t\in\Delta_T}\left\|r_{st}^i-\widetilde r_{st}^{\,i}\right\|\omega(s,t)^{-(N+1)\gamma}.\label{dZZ}
\end{eqnarray}
}
Now let
\[\Xi_{st}^i:=\langle \mathbf X_{st},\, L_i\big(\varphi_i(\mathbf Y)_s\big)\rangle \hbox{ and } \widetilde \Xi_{st}^{\,i}:=\langle \widetilde{\mathbf X}_{st},\, L_i\big(\widetilde{\varphi}_i(\widetilde{\mathbf Y})_s\big)\rangle,\]
and let $\Omega_{st}^i:=\Xi_{st}^i-\widetilde \Xi_{st}^{\,i}$.
From Lemma \ref{lem:deltaZ} we have
\begin{eqnarray}
\delta\Omega^i_{sut}&=&\left\langle\mathbf X_{ut},\, L_i\Big(R_{\mathbf X}\big(\varphi_i(\mathbf Y)\big)_{su}\Big)\right\rangle
-\left\langle\widetilde{\mathbf X}_{ut},\, L_i\Big(R_{\widetilde{\mathbf X}}\big(\widetilde{\varphi}_i(\widetilde{\mathbf Y})\big)_{su}\Big)\right\rangle\nonumber\\
&=&\sum_{j=1}^N\left\langle\pi_j^*(\mathbf X_{ut}),\, L_i\circ\pi_{j-1}\Big(R_{\mathbf X}\big(\varphi_i(\mathbf Y)\big)_{su}\Big)\right\rangle
-\sum_{j=1}^N\left\langle\pi_j^*(\widetilde{\mathbf X}_{ut}),\, L_i\circ\pi_{j-1}\Big(R_{\widetilde{\mathbf X}}\big(\widetilde{\varphi}_i(\widetilde{\mathbf Y})\big)_{su}\Big)\right\rangle\label{delta-Omega}
\end{eqnarray}
\begin{lemma}\label{est-delta-Omega}
The following estimate holds:
\begin{eqnarray}
\|\delta\Omega^i_{sut}\|&\le& \|L_i\|\bigg(\|\mathbf X\|_{\omega,\gamma}\,\|\varphi_i(\mathbf Y)-\widetilde\varphi_i(\mathbf Y)\|
+\|\mathbf X-\widetilde{\mathbf X}\|_{\omega,\gamma}\,\|\widetilde\varphi_i(\mathbf Y)\|\nonumber\\
&&\hskip 48mm +\|\widetilde{\mathbf X}\|_{\omega,\gamma}\,d\left(\widetilde\varphi_i(\mathbf Y),\, \widetilde\varphi_i(\widetilde{\mathbf Y})\right)\bigg)\sum_{j=1}^N\omega(u,t)^{j\gamma}\omega(s,u)^{(N+1-j)\gamma}.\nonumber
\end{eqnarray}
\end{lemma}
\begin{proof}
From \eqref{delta-Omega} we get
\begin{eqnarray*}
\delta\Omega^i_{sut}&=&\sum_{j=1}^N\bigg(\left\langle\mathbf X_{ut},\, L_i\Big(R_{\mathbf X}\big((\varphi_i-\widetilde\varphi_i)(\mathbf Y)\big)_{su}\Big)\right\rangle\\
&&\hskip 12mm + \left\langle\mathbf X_{ut}-\widetilde {\mathbf X}_{ut},\, L_i\Big(R_{\mathbf X}\big(\widetilde\varphi_i(\mathbf Y)\big)_{su}\Big)\right\rangle\\
&&\hskip 12mm + \left\langle\widetilde {\mathbf X}_{ut},\, L_i\Big(R_{\mathbf X}\big(\widetilde\varphi_i(\mathbf Y)\big)_{su}-R_{\widetilde{\mathbf X}}\big(\widetilde\varphi_i(\widetilde{\mathbf Y})\big)_{su}\Big)\right\rangle\bigg),
\end{eqnarray*}
from which Lemma \ref{est-delta-Omega} proceeds directly.
\end{proof}
\begin{corollary}\label{est-omega}
The following estimate holds:
\begin{equation}\label{est-omega-one}
\|r^i_{st}-\widetilde r\,^i_{st}\|\le NC \|L_i\|\bigg(\|\mathbf X\|_{\omega,\gamma}\,\|\varphi_i(\mathbf Y)-\widetilde\varphi_i(\mathbf Y)\|
+\|\mathbf X-\widetilde{\mathbf X}\|_{\omega,\gamma}\,\|\widetilde\varphi_i(\mathbf Y)\|+\|\widetilde{\mathbf X}\|_{\omega,\gamma}\,C_{\widetilde\varphi_i}(T)  d\left(\mathbf Y,\, \widetilde{\mathbf Y}\right)\bigg)\,\omega(s,t)^{(N+1)\gamma},
\end{equation}
where $C$ is a combinatorial constant.
\end{corollary}
\begin{proof}
It is proven by the same way than Theorem \ref{ixi-xi} and Corollary \ref{ixi-xi-bis}. It is a direct consequence of Theorem \ref{delta-xi} applied to $\Omega^i$, and Theorem \ref{thm:gen-sewing} applied to
\[\frac{\Omega^i}{\|L_i\|\bigg(\|\mathbf X\|_{\omega,\gamma}\,\|\varphi_i(\mathbf Y)-\widetilde\varphi_i(\mathbf Y)\|
+\|\mathbf X-\widetilde{\mathbf X}\|_{\omega,\gamma}\,\|\widetilde\varphi_i(\mathbf Y)\|+\|\widetilde{\mathbf X}\|_{\omega,\gamma}\,d\left(\widetilde\varphi_i(\mathbf Y),\, \widetilde\varphi_i(\widetilde{\mathbf Y})\right)\bigg)},\]
in the particular case when all controls $\omega_{1,j}$ and $\omega_{2,j}$ coincide with $\omega$. The factor $d\left(\widetilde\varphi_i(\mathbf Y),\, \widetilde\varphi_i(\widetilde{\mathbf Y})\right)$ in the denominator can be replaced by $C_{\widetilde\varphi_i}(T)  d\left(\mathbf Y,\, \widetilde{\mathbf Y}\right)$ thanks to Proposition \ref{d-varphi}.
\end{proof}
\noindent Plugging \eqref{est-omega-one} into \eqref{dZZ}, and using Proposition \ref{main-est-varphi} together, with \eqref{Kphi} we get
\begin{eqnarray}
d(\mathbf Z,\,\widetilde{\mathbf Z})&\le& \|\mathbf Y_0-\widetilde{\mathbf Y}_0\|+\omega(0,T)^\gamma\sum_{i=1}^d \|L_i\|\bigg(\left\|(\varphi_i-\widetilde\varphi_i)(\mathbf Y)\right\|+C_{\widetilde\varphi_i}(T)d\left(\mathbf Y,\,\widetilde{\mathbf Y}\right)\bigg)\nonumber\\
&& +NC\omega(0,T)^\gamma\sum_{i=1}^d \|L_i\|\bigg(\|\mathbf X\|_{\omega,\gamma}\,\|\varphi_i(\mathbf Y)-\widetilde\varphi_i(\mathbf Y)\|
+\|\mathbf X-\widetilde{\mathbf X}\|_{\omega,\gamma}\,\|\widetilde\varphi_i(\mathbf Y)\|+\|\widetilde{\mathbf X}\|_{\omega,\gamma}\,C_{\widetilde\varphi_i}(T)d\left(\mathbf Y,\, \widetilde{\mathbf Y}\right)\bigg).\nonumber\\
&\le& \|\mathbf Y_0-\widetilde{\mathbf Y}_0\|+\omega(0,T)^\gamma\left(NC\|\mathbf X\|_{\omega,\gamma}+1\right)\sum_{i=1}^d \|L_i\|\left\|(\varphi_i-\widetilde\varphi_i)(\mathbf Y)\right\|\nonumber\\
&& \hskip 19mm+NC\omega(0,T)^\gamma\|\mathbf X-\widetilde{\mathbf X}\|_{\omega,\gamma}\sum_{i=1}^d\|L_i\|\left\|\widetilde\varphi_i(\mathbf Y)\right\|\nonumber\\
&&\hskip 19mm +\omega(0,T)^\gamma\left(NC\|\widetilde{\mathbf X}\|_{\omega,\gamma}+1\right)\sum_{i=1}^d\|L_i\|C_{\widetilde\varphi_i}(T)d\left(\mathbf Y,\, \widetilde{\mathbf Y}\right)\nonumber\\
&\le& \|\mathbf Y_0-\widetilde{\mathbf Y}_0\|+\omega(0,T)^\gamma\left(NC\|\mathbf X\|_{\omega,\gamma}+1\right)\sum_{i=1}^d \|L_i\|C_{\varphi_i-\widetilde\varphi_i}(T)\|\mathbf Y\|\nonumber\\
&& \hskip 19mm+NC\omega(0,T)^\gamma\|\mathbf X-\widetilde{\mathbf X}\|_{\omega,\gamma}\sum_{i=1}^d\|L_i\|\left\|\widetilde\varphi_i(\mathbf Y)\right\|\nonumber\\
&&\hskip 19mm +\omega(0,T)^\gamma\left(NC\|\widetilde{\mathbf X}\|_{\omega,\gamma}+1\right)\sum_{i=1}^d\|L_i\|C_{\widetilde\varphi_i}(T) d\left(\mathbf Y,\, \widetilde{\mathbf Y}\right)\nonumber\\
&\le& \|\mathbf Y_0-\widetilde{\mathbf Y}_0\|+\omega(0,T)^\gamma\left(NC\|\mathbf X\|_{\omega,\gamma}+1\right)\sum_{i=1}^d \|L_i\|\Big(CK(T)^{N^2}\|\varphi_i-\widetilde\varphi_i\|_{K(Y,Y')}^N\sum_{i,j=0}^{N-1}\|\mathbf X\|_{\omega,\gamma}^i\|\mathbf Y\|^j\Big)\|\mathbf Y\|\nonumber\\
&& \hskip 19mm+NC\omega(0,T)^\gamma\|\mathbf X-\widetilde{\mathbf X}\|_{\omega,\gamma}\sum_{i=1}^d\|L_i\|\left\|\widetilde\varphi_i(\mathbf Y)\right\|\nonumber\\
&&\hskip 19mm +\omega(0,T)^\gamma\left(NC\|\widetilde{\mathbf X}\|_{\omega,\gamma}+1\right)\sum_{i=1}^d\|L_i\|C_{\widetilde\varphi_i} (T) d\left(\mathbf Y,\, \widetilde{\mathbf Y}\right).
\label{pre-ult}
\end{eqnarray}
\noindent We are now ready to state the Universal Limit Theorem:
\begin{theorem}\label{thm:ultfinal}
Let $\mathbf Z$ be the unique solution of the initial value problem \eqref{ivplifted} in $\mathcal V_{\mathbf X}^{[0,T']}$, and let $\widetilde{\mathbf Z}$ be the unique solution of the similar initial value problem in $\mathcal V_{\widetilde{\mathbf X}}^{[0,T']}$, with $\varphi_i$ replaced by $\widetilde{\varphi}_i$, where $T'\in[0,T]$ is small enough. Then we have
\[d(\mathbf Z,\, \widetilde{\mathbf Z})\le A\|\mathbf Z_0-\widetilde{\mathbf Z}_0\|+\sum_{i=1}^dB_i\|\varphi_i-\widetilde\varphi_i\|_{K(Z,Z')}^N
+E\|\mathbf X-\widetilde{\mathbf X}\|_{\omega,\gamma},\]
where $A$, $B_i$ and $E$ depend on $T',\mathbf Z, \mathbf X, \widetilde{\mathbf X}$ and the derivatives of the $\varphi_i$'s and $\widetilde{\varphi}_i$'s up to order $N$.
\end{theorem}
\begin{proof}
Set $\mathbf Z=\mathbf Y$ and $\widetilde{\mathbf Z}=\widetilde{\mathbf Y}$ in \eqref{pre-ult}. By replacing $T$ with $T'$ small enough such that
\[\omega(0,T')^\gamma\left(NC\|\widetilde{\mathbf X}\|_{\omega,\gamma}+1\right)\sum_{i=1}^d\|L_i\|C_{\widetilde\varphi_i}(T')< 1,\] we have
\ignore{
\[d(\mathbf Z,\, \widetilde{\mathbf Z})\le A\|\mathbf Z_0-\widetilde{\mathbf Z}_0\|+\sum_{i=1}^dB_i\|\varphi_i-\widetilde\varphi_i(\mathbf Z)\|+C\|\mathbf X-\widetilde{\mathbf X}\|_{\omega,\gamma}.\]
Finally applying Corollary \ref{main-est-varphi2} yields
}
\[d(\mathbf Z,\, \widetilde{\mathbf Z})\le A\|\mathbf Z_0-\widetilde{\mathbf Z}_0\|+\sum_{i=1}^dB_i\|\varphi_i-\widetilde\varphi_i\|_{K(Z,Z')}^N
+E\|\mathbf X-\widetilde{\mathbf X}\|_{\omega,\gamma}.\]
This completes the proof.
\end{proof}

\end{document}